\documentclass[a4paper,12pt,oneside,reqno]{amsart}
\usepackage{amsfonts,amssymb,amsmath,amsthm}
\usepackage[a4paper]{typearea}

\newtheorem{theorem}{Theorem}[section]
\newtheorem{lemma}[theorem]{Lemma}
\newtheorem{proposition}[theorem]{Proposition}
\newtheorem{corollary}[theorem]{Corollary}

\theoremstyle{definition}
\newtheorem{definition}[theorem]{Definition}

\newtheorem{condition}{Condition}
\newtheorem*{conditiona}{Condition}

\theoremstyle{remark}
\newtheorem*{remark}{Remark}

\numberwithin{equation}{section}

\def\Min{{\mathrm{min}}}
\def\dist{\mathop{\rm dist}\nolimits}

\def\bbone{{\mathchoice {\rm 1\mskip-4mu l} {\rm 1\mskip-4mu l}
          {\rm 1\mskip-4.5mu l} {\rm 1\mskip-5mu l}}}

\def\Asl{\mathop{\mathsf{Asl}}\nolimits}
\def\out{{\mathrm{out}}}
\def\d{{\mathrm{d}}}

\begin{document}

\title[ ]{The Arcsine law as a universal aging scheme for trap models}
\author[G. Ben Arous]{G\'erard Ben Arous}
\address{G\'erard Ben Arous\\ 
  Courant Institute of Mathematical Sciences\\
  New York University\\
  251 Mercer Street\\
  New York, N.Y. 10012-1185\\
  \and
  \'Ecole Polytechnique F\'ederale de Lausanne\\
  1015 Lausanne\\
  Swit\-zerland
}
\email{gerard.benarous@epfl.ch}

\author[J. \v Cern\'y]{Ji\v r\'\i~\v Cern\'y$^*$}
\address{Ji\v r\'\i~\v Cern\'y\\ 
  \'Ecole Polytechnique F\'ederale de Lausanne\\
  1015 Lausanne\\
  Switzerland
  }
\email{jiri.cerny@epfl.ch}
\thanks{$^*$ Work supported by the DFG Research Center ``\textsc{Matheon}''}

\subjclass[2000]{82D30, 82C41, 60F17}
\keywords{Aging, trap model, L\'evy process, random walk, time change}

\date{\today}
\begin{abstract}
  We give a general proof of aging for trap models using the arcsine law 
  for stable subordinators. This proof is based on 
  abstract conditions on the potential theory of the underlying 
  graph and on the randomness of the trapping landscape. We apply this 
  proof to aging for trap models on large two-dimensional tori and for 
  trap dynamics of the Random Energy Model on a broad range of 
  time scales. 
\end{abstract}

\maketitle

\section{Introduction} 
\label{s:introduction}

We establish in this article a general mechanism explaining the 
phenomenon of aging for some important dynamics in random media, 
i.e.~the ``trap models''. This general scheme is based on the classical 
arcsine law for stable subordinators. The general context is the 
following. Consider $Y(i)$ a discrete-time Markov chain on a discrete 
space $\mathcal V$ with transition kernel $p(x,y)$. Assume that a 
function $\boldsymbol \tau$ is given from the state space $\mathcal V$ 
to $(0,\infty)$. $\boldsymbol \tau$ will later be assumed to be random. It 
should be seen as a random landscape or a random scenery. Consider then 
the following sampling process: every time the Markov chain $Y$ is at 
$x\in \mathcal V$ 
it collects an (independent) exponential random variable with 
mean $\tau_x$. More precisely let $(e_i:i\ge 0)$ be an 
independent collection of i.i.d.~exponential mean-one random variables 
and define for any $k\ge 0$
\begin{equation}
  S(u)= \sum_{i=0}^{\lfloor u \rfloor -1} e_i \tau_{Y(i)}.
\end{equation}

The question we address first is the following: what is the behaviour 
of the process $S(u)$ for large $u$'s, when the landscape 
$\boldsymbol \tau$ is  highly heterogeneous, typically (but not 
  necessarily) when the $\tau_x$ are i.i.d.~and heavy tailed? We are 
not interested here in the case where $\tau$ is random but reasonably 
tame, which is usually studied under the name of Random Walk in 
Random Scenery (for recent important 
  results see \cite{AC06,GHK06,GKS06} and references therein). We will, 
on the contrary, isolate general conditions bearing both on the 
distribution of the random landscape $\boldsymbol \tau$ and on the 
potential theory of the chain $Y$ which will ensure that the process 
$S(u)$ can be approximated, in appropriate large time scales, by a 
stable subordinator. This convergence result will enable us to give a 
general mechanism explaining aging for the so-called Random Hopping 
Times dynamics or trap models, i.e.~for the continuous time Markov 
chain $X(t)$ whose jump rates are given by:
\begin{equation}
  w^{\boldsymbol \tau }_{xy}= \tau(x)^{-1} p(x,y).
\end{equation}
Indeed $X(t)$ is a time change of $Y(k)$ and the time change is given by 
(the inverse of) $S(u)$, which we call here the \textit{clock process}:
\begin{equation}
  X(t)=Y(k)\qquad \text{for all $t\in[S(k),S(k+1))$}.
\end{equation}

We will give a limit theorem stating that the clock process is close to 
a stable subordinator. We will then use the classical arcsine law, 
which gives the probability that a deterministic interval $(a,b)$ does not 
intersect the range of a subordinator as a function of the ratio $a/b$, 
in order to estimate the probability
\begin{equation}
  R(t_w,t_w+t;\boldsymbol \tau ):=
  \mathbb P[X(t_w+t)=X(t_w)|\boldsymbol \tau]
\end{equation}
for large times $t_w$ and $t$ in appropriate time scales. In 
particular, we will show that it is asymptotically a function of the 
ratio of these times $t_w/t$, which is what is usually called aging 
for the chain $X(t)$ in the statistical physics literature. 

The class of examples we do have in mind is the simple case where $Y$ 
is the standard random walk on a (connected) graph 
$G=(\mathcal V,\mathcal E)$, i.e.~where the transition probability is 
given by:
\begin{equation}
  p(x,y)= 1/d_x \qquad\text{if  $x$ and $y$ 
    are neighbours on the graph, $\langle x,y\rangle \in\mathcal  E$}
\end{equation}
and $d_x$ is the degree of the vertex $x$. This class of examples 
comes under the name of ``trap models''  in statistical mechanics of 
disordered media, and has been introduced by J.P Bouchaud (see 
  \cite{Bou92}). We refer to our lecture notes \cite{BC06c} for a 
complete survey and a more extensive bibliography. 

Let us summarise very briefly here the state of known results about 
aging for trap models. The trap models have been already 
studied on $\mathbb Z^d$ and on ``mean-field'' objects, that is on 
large complete graphs. Aging was first proved for large complete 
graphs in \cite{Bou92,BD95}. This case was seen by Bouchaud as a good 
ansatz for the dynamics of the simplest spin-glass,  the Random 
Energy Model (REM). It corresponds to the simplest case, where the 
Markov chain $Y(i)$ is simply a sequence of i.i.d. random variables  
uniformly distributed on a large finite set. Aging was then proved 
for the longest possible time scales for the REM dynamics in 
\cite{BBG03,BBG03b}, with a hard proof based on renewal theory. It was 
also proved  for the trap model on $\mathbb Z$ with a proof based on 
a direct scaling limit argument in \cite{FIN02} and \cite{BC05}. 
Finally, aging was proved  for the model on $\mathbb Z^d$, $d\ge 2$, 
on the shortest possible time scale \cite{BCM06,Cer03}. The  proof 
there is based on a difficult coarse-graining procedure. 

The striking fact is that these aging results are identical for 
$\mathbb Z^d$, $d\ge 2$ and the large complete graph, or the REM. In 
other terms, the mean-field results are valid from infinite dimension 
down to dimension $2$. 

The new approach we give here is based on what 
we have learnt from these examples and  has the following advantages

1. It shows very clearly how aging is based on the arcsine law, 
isolating the natural interplay between the potential theory of the 
chain and the randomness of the landscape. It also shows that even 
though Bouchaud's ansatz (the model on the complete graph) is not 
universal in finite dimensions, one of its features is, at least for 
$d\ge 2$, i.e.~the nature of the clock process. The two-time function $R$ 
being insensitive to anything but the range of the clock process, this is enough 
to imply aging, through the arcsine law, as soon as the approximation 
by a stable subordinator is valid.

2. It allows us to give aging results in a broad range of time scales. 
For instance on very long time scales in finite dimensions ($d\ge 2$), or 
in short time scales for the REM. We will see that there is a natural 
range of space scales (i.e.~level sets of $\tau$) and time scales in 
which the scheme based on the arcsine law applies. We will also see 
that it is possible in our scheme to have varying exponents in the 
arcsine law in varying time scales for the same model (the REM for 
  instance). Moreover this will show that the aging phenomenon is a 
question about the \textit{transient part of relaxation} to equilibrium and 
not necessarily related to equilibrium questions. The most striking 
illustration being the fact that we prove that aging can occur for 
the REM \textit{above the critical temperature} (where equilibrium questions are 
  trivial).

3. It allows us to think of random landscapes which would not be 
i.i.d.,~and thus open the possibility of studying aging of trap 
dynamics of more relevant spin glasses than the REM. At least with a good 
dose of optimism, and naturally for time scales short enough to not 
yet feel the model-specific equilibrium features.

\medskip

Let us now be more specific, and describe precisely the questions we 
address here. Even though everything we prove in this paper is valid 
in the general context described above with no change 
at all, we will restrict the exposition to the case of trap models 
for the sake of simplicity. We will nevertheless expand a bit the 
framework by considering a sequence of graphs rather than one fixed 
graph (in order to accommodate also the large complete graphs of 
  Bouchaud's ansatz or the hypercube in high dimensions needed for 
  the REM). Consider thus a sequence of connected graphs  
$G_n=(\mathcal V_n, \mathcal E_n)$, $n\in \mathbb N$, with the vertex 
set $\mathcal V_n$ and the edge set $\mathcal E_n$. 
Let $(Y_n(j):j\ge 0)$ be a discrete-time simple random walk 
on $G_n$.
For each vertex $x\in \mathcal V_n$, let $\tau_x$ 
be a non-negative real number, which we call the depth of the trap at 
$x$, and denote by $\boldsymbol \tau_n$ the collection of these 
depths, $\boldsymbol \tau_n =\{\tau_x:x\in \mathcal V_n\}$. We will 
assume that $\boldsymbol \tau_n$ is a sequence of ``random 
environments'', i.e.~random variables with distribution $\mu_n$ on 
$[0,\infty)^{\mathcal V_n}$. We suppose that $\mu_n$ are defined on a 
common probability space, so that we can consider a.s. convergence. 
Note that we do not assume a priori that the $\tau_x$'s are 
i.i.d.,~even though in the classical examples the $\tau_x$'s are 
i.i.d.~and heavy tailed \cite{BC06c}.

Given the environment $\boldsymbol \tau_n$, we define the trap 
model as a continuous-time  Markov process $X_n(\cdot)$ with state 
space $\mathcal V_n$ whose transition rates are given by 
\begin{equation}
  \label{e:dynamics}
  w_{xy}^{\boldsymbol \tau }=
  \begin{cases}
    (d_x \tau_x)^{-1}&\text{if $\langle x,y\rangle \in \mathcal E_n$,}\\
    0&\text{otherwise}.
  \end{cases}
\end{equation}
Here $d_x$ stands for  the degree of $x$ in the graph $G_n$. In 
words, $X_n$ waits at $x$ an exponentially distributed time with mean 
$\tau_x$ and then it jumps to one of the neighbours of $x$ with the equal 
probability $d_x^{-1}$.  We write $\mathbb P_x$ for the distribution 
of $X_n$ conditioned on $X_n(0)=x$. Usually, we will consider $X_n$ 
to be started at some arbitrary fixed vertex that does not depend on 
$\boldsymbol \tau_n $. This vertex is  denoted by $\boldsymbol 0$, we 
write $\mathbb P=\mathbb P_{\boldsymbol 0}$.

As we have already explained, $X_n$ is a time change of $Y_n$. 
Indeed,  define the {\em clock process} $S_n(u)$, $u\ge 0$ by
\begin{equation}
  \label{e:Sdef}
  S_n(u):=\sum_{i=0}^{\lfloor u\rfloor-1}e_i \tau_{Y_n(i)}.
\end{equation}
So that,
\begin{equation}
  X_n(t)=Y_n(j)\qquad \text{for all $t\in[S_n(j),S_n(j+1))$}.
\end{equation}

To study aging we need choose a two-time function that 
reflects the behaviour of the system in the time interval 
$[t_w, t_w+t]$. The most natural two-time function for the trap models is the 
probability that at both times $t_w$ and $t_w+t$ the system is in the 
same state,
\begin{equation}
  R_n(t_w,t_w+t;\boldsymbol \tau_n ):=
  \mathbb P[X_n(t_w+t)=X(t_w)|\boldsymbol \tau_n].
\end{equation}
There are other possible choices for the two-time function (see 
  \cite{BC06c}). We will however not consider them here.

\begin{definition}
  We say that the function $R_n$ exhibits aging  if for some sequence 
  $t(n)$ such that $\lim_{n\to\infty}t(n)=\infty$ it satisfies
  \begin{equation}
    \label{e:defaging}
    \lim_{n\to\infty}R_n(t(n),(1+\theta )t(n);\boldsymbol \tau_n 
      )=R(\theta )
  \end{equation}
  for all $\theta >0$ and some non-trivial function $R(\theta )$. We 
  call $R(\theta )$ the aging function. 
\end{definition}

As we have already remarked, in all  cases where aging of $R_n$ was proved,  
with the exception of the one-dimensional case, the limiting function  
$R(\theta )$ was given by the arcsine law for L\'evy processes%
\footnote{In \cite{BBG03b} a more complicated limiting procedure than in 
  \eqref{e:defaging} was used and the two-time function $R$ was 
  slightly different, we will make more comments about this issue later.},
\begin{equation}
  \label{e:arcsinuslaw}
  R(\theta )=\Asl_\alpha (1/1+\theta ),
\end{equation}
where $\Asl_\alpha (u)$ stands for the distribution function of the 
generalised arcsine law with parameter $\alpha $,
\begin{equation}
  \Asl_\alpha (z):=
 \frac{\sin \alpha \pi }{\pi } \int_0^{z}u^{\alpha 
	-1}(1-u)^{-\alpha }\,\d u.
\end{equation}
Note that $\Asl_\alpha (a/b)$ is equal to the probability that the 
range of an $\alpha $-stable subordinator does not intersect the interval 
$[a,b]$ \cite{Ber96}. 

The aim of this paper is to give a set of possibly simple conditions 
that guarantee for general graphs $G_n$ and time scales $t(n)$ the same 
behaviour, that is the convergence of $R_n(t(n),(1+\theta )t(n))$ 
to $\Asl_\alpha(1/1+\theta ) $.

We will give first a set of four general conditions (A)--(D) which 
ensure the convergence of the rescaled clock process to a stable 
subordinator. These four conditions are true for every known example 
(except, naturally, for the model on $\mathbb Z$ where the clock 
  process converge to a Kesten-Spitzer process \cite{KS79}, see also 
  \cite{BC06c}). We will then, in Section~\ref{s:mechanism}, give a 
set of four more general conditions 1--4 which are weaker but 
sufficient to ensure the convergence of the range of the clock 
process, which is enough to apply the arcsine law. In order to prove 
aging for the two-time function $R_n$ using this arcsine law we still 
need to impose two extra technical conditions, either on top of 
(A)--(D) or on top of 1--4.

\medskip

Let us introduce some notations useful in order to state our conditions.
Let $T_n$ be a stopping time for $Y_n$. We use $G_{T_n}^n(x,y)$, 
$x,y\in \mathcal V_n$, to denote the Green's function of $Y_n$, that 
is the mean time that $Y_n$ spends in $y$ before $T_n$ when started 
at $x$, 
\begin{equation}
  G^n_{T_n}(x,y)=\mathbb E_x\sum_{i=0}^{T_n-1}
  \mathbb \bbone\{Y_n(i)=y\}.
\end{equation}
For a set $A\subset \mathcal V_n$ we define its hitting time 
$H_n(A)$ as 
\begin{equation}
  H_{n}(A):=\inf\{i\ge 0:Y_n(i)\in A\}.
\end{equation}
To simplify the notation we write $G_A^n(\cdot,\cdot)$ for 
$G^n_{H(A)}(\cdot,\cdot)$. We define another two-time function 
$R_A^n(t_w,t_w+t)$ as the probability that $X$ does not visits any 
``fresh'' site in  $A$ during the observation interval: let $\ell_A(t_w)$ 
be the last time when $X_n$ visited $A$ before $t_w$,
\begin{equation}
  \ell_A(t_w)=\max\{t\le t_w:X_n(t)\in A\},
\end{equation}
and let 
$T=\inf\big\{s\ge t_w: X_n(s)\in A\setminus\{X_n(\ell_A(t_w))\}\big\}$, 
then 
\begin{equation}
  R_A^n(t_w,t_w+t;\boldsymbol \tau ) = 
  \mathbb P[T>t_w+t |\boldsymbol \tau ].
\end{equation}
A two-point function of this type was considered in \cite{BBG03}.

We further say that the set $A\subset \mathcal V_n$ is a 
\textit{Poisson cloud} on $\mathcal V_n$ with density $\rho \in (0,1)$ 
if  each site $x\in \mathcal V_n$ is in  $A$ with probability $\rho $ 
independently of all others, i.e.~if $A$ is a site-percolation 
process on $\mathcal V_n$. 

We can now formulate the first set of conditions that implies aging 
on the time scale $t(n)$. First, we need to 
control the behaviour of the random environment. 
\begin{conditiona}[\textbf{A}]
  For all $n$ the random environment $\boldsymbol \tau_n$ is i.i.d. 
  Further, there exist a depth scale $g(n)$, a density scale  $\rho (n)$  
  and a constant $\alpha \in (0,1)$ 
  such that $g(n)\to\infty$ and $\rho (n)\to 0$ as $n\to\infty$,
  and
  \begin{equation}
    \rho (n)^{-1}\mu_n[\tau_x\ge u g(n)]
    \xrightarrow{n\to\infty} u^{-\alpha },
  \end{equation}
  uniformly on all compact subsets of $(0,\infty)$.
  Moreover, there exist a constant $C$ such that for all $u>0$ and 
  $n\in \mathbb N$ 
  \begin{equation}
    \mu_n[\tau_x\ge u g(n)]\le C u^{-\alpha }\rho (n).
  \end{equation}
\end{conditiona}

The next two conditions control the motion of the simple random 
walk $Y_n$ between points of a Poisson cloud.

\begin{conditiona}[\textbf{B}]
    Let $A_n$ be a sequence of Poisson clouds on $\mathcal V_n$ 
    with densities $\rho  \rho (n)$, $\rho \in(0,\infty)$. 
    Then there exists a 
    constant $\mathcal K_G\in (0,\infty)$ independent of $\rho $ such 
    that  for a.e.~sequence $A_n$ 
    \begin{equation}
      \max_{x\in A_n}\Big|
      f(n)^{-1}G^n_{A_n\setminus\{x\}}(x,x)-
      \mathcal K_G
      \Big|
      \xrightarrow{n\to\infty}0,
    \end{equation}
    where the scale $f(n)$ is given by $f(n)=t(n)/g(n)$.
\end{conditiona}
\begin{conditiona}[\textbf{C}]
   There exists $\mathcal K_r\in (0,\infty)$ such that for all 
    $s>0$ and a.e.~sequence $A_n$ of Poisson clouds as in Condition (B)
    \begin{equation}
      \max_{x\in A_n\cup\{\boldsymbol 0\}}\bigg|
      \mathbb E_x\Big[\exp\Big(-\frac s{r(n)} H_n(A_n\setminus\{x\})\Big)\Big]-
      \frac{\mathcal K_r\rho }{s+\mathcal K_r\rho }
      \bigg|
      \xrightarrow{n\to\infty}0,
    \end{equation}
    where $r(n) \rho (n)=f(n)$. In other words, 
    $H_n(A_n\setminus\{x\})/r(n)$ is asymptotically exponentially 
    distributed with mean $1/\mathcal K_r \rho $.
    (The scale $r(n)$ represents  the number of steps that $X_n$ 
    makes before time $t(n)$.) 
\end{conditiona}

Finally, we need one technical condition

\begin{conditiona}[\textbf{D}]
    There exists a large constant 
    $\mathcal K_s$ such that for  all $m>0$ and $n$ large
    \begin{equation}
      \sum_{x\in \mathcal V_n}\big(e^{\lambda_n G^n_{mr(n)}}-1\big)\le 
      \mathcal K_s \lambda_n \sum_{x\in \mathcal V_n}G^n_{mr(n)}(0,x)=
      \mathcal K_s \lambda_n m r(n),
    \end{equation}
    and $\sum_{n=1}^\infty \exp\big(-c \lambda_n f(n)\big)$ is 
    finite for all $c>0$.
\end{conditiona}

\begin{theorem}
  \label{t:pot}
  Assume that Conditions (A)--(D) hold. Then 

  (i) for a.e.~random environment $\boldsymbol \tau $ the 
  rescaled clock process ${t(n)}^{-1}S( r(n) \cdot )$ 
  converges to an $\alpha $-stable 
  subordinator weakly in the Skorokhod topology on $D([0,T])$ for all 
  $T>0$. 
  
  (ii) Further, define the set of deep traps
  $T_\varepsilon^M(n):=\{x\in \mathcal V_n:\varepsilon g(n)\le \tau_x< Mg(n)\}$. 
  Then $R^n_{T_\varepsilon^M(n)}$ ages: 
  \begin{equation}
    \lim_{\substack{\varepsilon \to0\\M\to \infty}}\lim_{n\to\infty}
    R^n_{T_\varepsilon^M(n)}(t(n),(1+\theta )t(n))=\Asl_\alpha 
    (1/1+\theta ).
  \end{equation}
\end{theorem}

\section{Graph-independent mechanism of proof}
\label{s:mechanism}

In this section we state the second set of the conditions. This set 
is  adjusted  to prove aging for 
the two-time function $R$ and its convergence to $\Asl_\alpha $. 
The conditions are more complicated 
than the Conditions (A)--(D) of the first set, and can be rather regarded 
as parts of a mechanism of a proof of aging. 

To understand these conditions it is useful to keep in mind the 
analogy with the sum of i.i.d.~non-negative $\alpha $-stable random 
variables. It is 
known fact that, after a proper renormalisation, this sum converge to 
an $\alpha $-stable subordinator. Moreover, the sum is typically 
dominated by a finite number of large contributions whose size 
depend on the number of terms.
We want to prove that the same holds for the clock process.

To formulate the second set of  conditions it is necessary to choose several 
objects that depend on the particular sequence $G_n$ and on the 
observation time scale $t(n)$. 

First, it is necessary to fix a (random) time $\xi_n$  up to which we 
observe $Y_n$. This time will serve  as an upper time scale up to 
which we observe $Y_n$. It must therefore be chosen large enough to 
ensure that (with high probability) $S_n(\xi_n)$ is larger than 
$(1+\theta )t(n)$ (see Condition~\ref{c:large} below). On the other 
hand, $\xi_n$ should be as small as possible to simplify the 
verification of the  other conditions. 

Second, a scale $g(n)$ for deep traps should be 
chosen according to $G_n$ and $t(n)$. This scale defines the set of 
the \textit{deep traps} by
\begin{equation}
  \label{e:deftop}
  T_\varepsilon^M(n):=\{x\in \mathcal V_n: \varepsilon g(n)\le 
    \tau_x < M g(n)\}.
\end{equation}
This set will 
determine the behaviour of the clock process $S_n$  at the  time scale 
$t(n)$. That is the clock process 
should be dominated by a very small number of relatively 
large contributions  due to visits of deep traps in 
$T_\varepsilon^M(n)$.  

We use 
$T_M(n):=T_M^\infty(n)=\{x\in \mathcal V_n: M g(n)\le \tau_x \}$ to 
denote the set of {\em very deep traps}. Similarly, we write 
$T^\varepsilon (n):=T^\varepsilon_0(n)=\{x\in \mathcal V_n:  \tau_x < \varepsilon  g(n)\}$ 
for the set of {\em shallow traps}. To justify the analogy with the 
sum of i.i.d.~random variables, the contribution of these two sets 
should be negligible. This is the content of the first two 
conditions. Condition~\ref{c:shallow}  states that shallow traps are 
irrelevant because the time spent in them is too small. 
Condition~\ref{c:deep} states that very deep traps are irrelevant 
because they are not seen by the Markov chain. 
\begin{condition}
  \label{c:shallow}
  There is a function $h(\varepsilon )$  satisfying 
  $\lim_{\varepsilon \to 0} h(\varepsilon )=0$, such that for 
  a.e.~realisation of $\boldsymbol \tau := \{\boldsymbol \tau_n:n>0\}$  
  and for all $n$ large enough
  \begin{equation}
    \mathbb E\Big[\sum_{i=0}^{\xi_n}e_i \tau_{Y_n(i)}\bbone\{Y_n(i)\in 
	T^\varepsilon (n)\}\Big|\boldsymbol \tau \Big]\le h(\varepsilon )t(n).
  \end{equation}
  That is,  the expected time spent in shallow traps before $\xi_n$ is 
  small with respect to $t(n)$.
\end{condition}

The second conditions ensures the negligibility of the very deep 
traps. Recall that $H_{n}(A)$  
denotes the hitting time of the set $A$ by  $Y_n$.
\begin{condition}
  \label{c:deep}
  Given $\xi_n$, for any $\delta >0$ there exists $M$ large enough such 
  that for a.e.~realisation of $\boldsymbol \tau $ and for all $n$ large 
  \begin{equation}
    \mathbb P\big[H_n(T_M(n))\le \xi_n\big|\boldsymbol \tau \big]\le \delta.
  \end{equation}
\end{condition}

We need other definitions to state conditions that  
guarantee the existence of the limit in \eqref{e:defaging}. 
First, let $r_n(j)$ be the sequence of times when a new deep trap is 
visited, $r_n(0)=0$, and 
\begin{equation}
  r_n(i)=\min\big\{j>r_n(i-1): Y_n(j)\in T_\varepsilon^M(n)
    \setminus\{Y_n(r_n(i-1))\}\big\}. 
\end{equation}
We use $\zeta_n$ to denote the largest $j$ such that $r_n(j)\le \xi_n$,
\begin{equation}
  \label{e:defzeta}
  \zeta_n:=\max \{j:r_n(j)\le \xi_n\}.
\end{equation} 
We define the process $U_n(j)$ that records the trajectory of $Y_n$ (and 
  thus of~$X_n$) restricted to the deep traps, 
\begin{equation}
  \label{e:Undef}
  U_n(j):=Y_n(r_n(j)),\qquad j\in \mathbb N_0.
\end{equation}
Finally, let $s_n(j)$ be the time that $X_n$ spends at site $U_n(j)$ 
between steps $r_n(j)$ and $r_n(j+1)$, 
\begin{equation}
  \label{e:defsn}
  s_n(j):=\sum_{i=r_n(j)}^{r_n(j+1)}e_i \tau_{Y_n(i)} 
  \bbone\{Y_n(i)=U_n(j)\},\qquad
  \text{$j<\zeta_n$}.
\end{equation}
It is easy to observe that $s_n(j)$ has an exponential distribution 
with mean 
\begin{equation}
  \tau_{U_n(j)}G^n_{T_\varepsilon^M\setminus\{U_n(j)\}}(U_n(j),U_n(j)).
\end{equation}

Since Conditions \ref{c:shallow} and \ref{c:deep} ensure that the 
visits of the deep traps determine the behaviour of the time change $S_n(j)$, 
the sum $\sum_{i=1}^{j-1} s_n(i)$ can be considered as a good 
approximation of $S_n(r_n(j))$.  We would like to show that the $s_n(j)$ 
become independent as $n\to \infty$, and that they have an 
appropriate tail behaviour. To this end, we define 
$(\sigma_\varepsilon^M(i),i\in \mathbb N)$ as a sequence of 
i.i.d.~random variables taking values between $\varepsilon $ and $M$ 
with common distribution function
\begin{equation}
  \label{e:defsigman}
  \mathbb P[\sigma_\varepsilon^M(i)\le u]=
  \frac {\varepsilon^{-\alpha }-u^{-\alpha }}
  {\varepsilon^{-\alpha }-M^{-\alpha }}=:
  \frac{\varepsilon^{-\alpha }-u^{-\alpha }} {p_\varepsilon^M }
  , \qquad u\in [\varepsilon ,M].
\end{equation}
Let $(\hat e_i,i\in \mathbb N)$ be a sequence of mean-one i.i.d.~exponential 
random variables that are independent of 
$\sigma_\varepsilon^M$, and let 
$s_\infty(i):=\hat e_i \sigma_\varepsilon^M(i)$. For notational 
convenience we define $s_n(j)=s_\infty(j)$ for all $j\ge \zeta_n$.

The following conditions will ensure that the limit in 
\eqref{e:defaging} is given by the arcsine law \eqref{e:arcsinuslaw}.
First, we need the asymptotic independence and the proper tail behaviour:
\begin{condition}
  \label{c:conv}
  There exists a constant $\mathcal K>0$ such that
  for all $\varepsilon$, $M$ and 
  for a.e.~$\boldsymbol \tau $,  the sequence 
  $(s_n(j)/t(n),j\in \mathbb N)$ converges as $n\to \infty$ in law to the sequence 
  of i.i.d.~random variables $(\mathcal K s_\infty(j),j\in \mathbb N)$. 
\end{condition}
We need also to ensure that $S_n(r_n(\zeta_n))$ is 
larger than $(1+\theta)t(n)$ with a large probability. Since 
$r_n(\zeta_n)\ge \sum_{i=1}^{\zeta_n-1}s_n(i)$, and $s_n(i)$ are 
easier to control than $S_n(r_n(j))$ we require
\begin{condition}
  \label{c:large}
  For a.e.~$\boldsymbol \tau $ and for any fixed 
  $\theta >0$, $\delta >0$ it is possible to choose $\xi_n$ such that for all 
  $\varepsilon $ small and $M$ large enough, and for $\zeta_n$ 
  defined in \eqref{e:defzeta}
  \begin{equation}
    \label{e:c:large}
    \mathbb P\Big[\sum_{i=1}^{\zeta_n-1}s_n(i)\ge (1+\theta 
	)t(n)\Big|\boldsymbol \tau \Big]\ge 1-\delta .
  \end{equation}
\end{condition}

\medskip

The next pair of conditions is, in principle, necessary only for a 
``post-processing''.  If they are not verified, it is possible to  
prove aging for the {\em $T_\varepsilon^M(n)$-dependent} two-time 
function $R^n_{T_\varepsilon^M(n)}$. Observe that  this function can 
be also written as
\begin{equation}
  \begin{split}
    R^n_{T_\varepsilon^M(n)}(t_w,t_w+t;\boldsymbol \tau )&=
    \mathbb P\big[\exists j :
      S_n(r_n(j))\le t_w < t_w +t< S_n(r_n(j+1))\big|\boldsymbol \tau \big]\\
    &=
    \mathbb P [ \{S_n(j):j\in \mathbb N\}\cap (t_w,t_w+t]=\emptyset
      |\boldsymbol \tau ].
  \end{split}
\end{equation}

To prove aging for the two-point function $R$  we need to know that for 
any time $t'$ between $S_n(r_n(j))$ and $S_n(r_n(j+1))$ the probability that 
$X_n(t')=U_n(j)$ is large. 
For a formal statement of this claim we need some more definitions.
Let $t'_n$ be a deterministic time sequence satisfying 
$t(n)/2\le t'_n\le (1+\theta )t(n)$, and let $\delta >0$.
We define $j_n\in \mathbb N$ by 
\begin{equation}
  \label{e:barjn}
  S_n(r_n(j_n))\le t'_n\le S_n(r_n(j_n+1))-\delta t(n), 
\end{equation}
and $j_n=\infty$ if \eqref{e:barjn} is not satisfied for any integer.
Let $A_n(\delta )$ be the event 
\begin{equation}
  \label{e:defAn}
  A_n(\delta ):=\{ 0< j_n<\zeta_n\}.
\end{equation}
\begin{condition}
  \label{c:post}
  For any $\delta$ it is possible to choose $\varepsilon $ small  and $M$ large 
  enough such that for a.e.~$\boldsymbol \tau $ and all $n$ large enough
  \begin{equation}
  \mathbb P[X_n(t'_n)=U_n(j_n)|A_n(\delta ),\boldsymbol \tau ]\ge 1-\delta.
  \end{equation}
\end{condition}

The last condition that we need to prove aging for $R$ excludes 
repetitions in the sequence $U_n$.
\begin{condition}
  \label{c:noreturn}
  \label{c:last}
  For any fixed $\varepsilon$ and $M$ and a.e.~$\boldsymbol \tau $
  \begin{equation}
    \lim_{n\to \infty}\mathbb P[\exists 0< i,j\le \zeta_n \text{ such 
	that } i\neq j \text { and } U_n(i)=U_n(j)|\boldsymbol\tau  ]=0.
  \end{equation}
\end{condition}

We now show how to use these six conditions  to prove the aging 
behaviour for the two-time functions  $R_n$, $R^n_{T_\varepsilon^M(n)}$.

\begin{theorem}
  \label{t:main}
  (i) Assume that Conditions \ref{c:shallow}--\ref{c:last} are 
  satisfied. Then for a.e.~realisation of the random environment 
  $\boldsymbol \tau $
  \begin{equation}
    \lim_{n\to\infty}R_n(t(n),(1+\theta )t(n);\boldsymbol \tau 
      )=\Asl_\alpha (1/1+\theta ).
  \end{equation}
  (ii) If only Conditions \ref{c:shallow}--\ref{c:large} hold, then the 
  same is valid for the two-time function $R^n_{T_\varepsilon^M(n)}$,
  \begin{equation}
    \lim_{\substack{\varepsilon \to0\\M\to \infty}}
    \lim_{n\to\infty}R^n_{T_\varepsilon^M(n)}(t(n),(1+\theta )t(n);\boldsymbol \tau 
      )=\Asl_\alpha (1/1+\theta ).
  \end{equation}
\end{theorem}


Before proving Theorems~\ref{t:pot} and~\ref{t:main} let us explain 
how both sets of conditions are related.
\begin{proposition}
  \label{p:rel}
  Conditions (A)--(D) imply Conditions~\ref{c:shallow}--\ref{c:large} 
  for the same scale $g(n)$ and for $\xi_n=m r(n)$ with some large 
  $\theta $-dependent constant $m$.  
\end{proposition}

We will use this proposition for both examples that we study later. 
That is to prove aging for $R_n$ we will always verify Conditions 
(A)--(D), \ref{c:post} and~\ref{c:noreturn}. 

Observe also that neither in 
Conditions~\ref{c:shallow}--\ref{c:last}, neither in 
Theorem~\ref{t:main} we  suppose that the $\tau_x$'s are i.i.d. 
This assumption is however used twice when proving 
Proposition~\ref{p:rel} as we will see later.
First, we will use the independence to verify Condition~\ref{c:shallow} from 
Condition~(D)
(see formulas \eqref{e:sh1}, \eqref{e:sh1a} below), that is to prove 
that the time spent in the shallow traps is small. We do however believe 
that Condition~\ref{c:shallow} stays valid  also for some dependent random 
environment. 
The second use of the independence is more substantial. It 
implies that the   geometrical 
structure of the set of the deep traps is particularly simple: it is 
a Poisson cloud. It is therefore easy to control, e.g., the minimal 
distance between deep traps or the relative size of the slices 
$T_\varepsilon^u$, $T_\varepsilon^M$. This control can be 
problematic when the $\tau_x$'s are dependent.   

Remark also that  Theorem~\ref{t:pot}(ii) is 
a simple consequence of Proposition~\ref{p:rel} and Theorem~\ref{t:main}(b). 
Therefore, we first prove Theorem~\ref{t:main}, then we verify 
Proposition~\ref{p:rel}. 
In the end we show   Theorem~\ref{t:pot}(i), that is the 
convergence of the clock process. Note that 
Conditions~\ref{c:shallow}--\ref{c:last} are not strong enough to 
imply directly such a convergence. We are however not aware of any particular 
case where Conditions~\ref{c:shallow}--\ref{c:last} hold and this 
convergence does not take place.

\begin{proof}[Proof of Theorem~\ref{t:main}]
  Let us define 
  \begin{equation}
    \tilde S(j)=\frac 1 {t(n)}\sum_{i=1}^{j-1} s_n(i)
  \end{equation}
  and let $E=E(n)$ be the event whose probability we are trying to estimate,
  \begin{equation}
    \label{e:defE}
    E(n)=\{X_n(t(n))=X_n((1+\theta )t(n))\}.
  \end{equation}
  We first explain the strategy of the proof. The most important 
  observation is that the clock process $S_n$ contains enough information 
  about $X_n$  to prove aging.
  Between times $S_n(r_n(j))$ and 
  $S_n(r_n(j+1))$ the process $X_n$ visits (possibly many times) 
  only one deep trap, $U_n(j)$, and it also visits many shallow 
  traps. Condition~\ref{c:post} ensures that if we pick a time $t$ 
  between $S_n(r_n(j))$ and $S_n(r_n(j+1))$, then $X_n(t)=U_n(j)$ 
  with a high probability. That means that if 
  $S_n(r_n(j))\le t(n)\le (1+\theta )t(n)<S_n(r_n(j+1))$ for some 
  $j\ge 1$, then $E$ holds with a  probability close to $1$. On 
  the other hand, if there is a  $j$ such that 
  $t(n)<S_n(r_n(j))<(1+\theta )t(n)$, then, using 
  Condition~\ref{c:noreturn}, $E$ can happen only if at both times $t(n)$ 
  and $(1+\theta )t(n)$ the process is in the same shallow trap. However, 
  this event has (again by Condition~\ref{c:post}) a very small 
  probability. Therefore, it is important to estimate the 
  probability that there is no $S_n(r_n(j))$ in the time interval of 
  interest. To this end, we will show that $\tilde S_n(j)$ is a 
  good approximation of $S_n(r_n(j))/t(n)$, and then we will estimate 
  the probability that there is no $\tilde S_n(j)$ in $[1,1+\theta ]$. 
  At the end of the proof we use these results to give a rigorous 
  version of the reasoning in this paragraph. 

  \medskip

  First, let us show that $\tilde S_n(j)$ approximates well 
  $S_n(r_n(j))/t(n)$, at least for all relevant indices $j\le\zeta_n$. 
  \begin{lemma}
    \label{l:approx}
    For all $\theta $ and $\delta >0$ there exist $\xi_n$, 
    $\varepsilon $, and $M$ such that $\boldsymbol \tau $-a.s.
    \begin{equation}
    \label{e:approx1}
      \mathbb P\Big[
	\max\Big\{\Big|\frac {S_n(r_n(j))}{t(n)}-\tilde S_n(j)\Big|:
	  S_n(r_n(j))\le (1+\theta )t(n)\Big\}>\delta\Big |
	\boldsymbol \tau  \Big]<\delta.
    \end{equation}
  \end{lemma}
  \begin{proof}
    By Condition~\ref{c:large} we can choose $\xi_n$ not depending on 
    $\varepsilon $ and $M$ such that 
    \begin{equation}
      \label{e:approx1a}
      \mathbb P[\tilde S_n(\zeta_n )\le 1+\theta ]\le \delta /2.
    \end{equation}
    Observing that $S_n(r_n(j))/t(n)-\tilde S_n(j)$ is positive and 
    increasing in $j$, it is sufficient to estimate 
    $S_n(r_n(\zeta_n))/t(n)-\tilde S_n(\zeta_n)$. However,
    \begin{equation}
      \frac{S_n(r_n(\zeta_n))}{t(n)}-\tilde S_n(\zeta_n)\le
      \sum_{i=0}^{\xi_n}
      \frac{e_i \tau_{Y_n(i)}}{t(n)}
      \bbone\{Y_n(i)\in T^\varepsilon (n)\cup T_M(n)\}.
    \end{equation}
  The contribution coming from $T^\varepsilon (n)$ can be bounded using 
  Condition~\ref{c:shallow} and Chebyshev inequality,
  \begin{equation}
    \label{e:approx2}
    \mathbb P\Big[t(n)^{-1} \sum_{i=0}^{\xi_n} e_i \tau_{Y_n(i)}
      \bbone\{Y_n(i)\in T^\varepsilon (n)\}\ge \delta/2\Big|\boldsymbol 
	\tau \Big]\le c 
      \delta^{-1}h(\varepsilon ).
  \end{equation}
  Since $h(\varepsilon )\to 0$ as $\varepsilon \to 0$, we can fix 
  $\varepsilon $ such that the last expression is bounded by 
  $\delta /4$. Similarly, using Condition~\ref{c:deep}, we can choose 
  $M$ such that 
  \begin{equation}
    \label{e:approx3}
    \mathbb P\Big[ \sum_{i=0}^{\xi_n} e_i \tau_{Y_n(i)}
      \bbone\{Y_n(i)\in T_M (n)\}\neq 0\Big | \boldsymbol \tau \Big]\le \delta /4.
  \end{equation}
  The lemma then follows combining \eqref{e:approx1a}, 
  \eqref{e:approx2} and \eqref{e:approx3}.
  \end{proof}

  We further compute the probability that an interval does not 
  contain any of the $\tilde S_n(j)$'s. 
  \begin{lemma}
    \label{l:levyconv}
    For all $0<a<b$ and for a.e.~$\boldsymbol \tau $
    \begin{equation}
      \lim_{\substack {\varepsilon \to 0 \\ M\to \infty}}\lim_{n\to \infty}
      \mathbb P\big[[a,b]\cap \{\tilde S_n(j), j\in \mathbb N\}=\emptyset\big|
	\boldsymbol \tau \big]=
      \Asl_\alpha (a/b).
    \end{equation}
  \end{lemma}
  \begin{proof}
    Let  $y_0=0$ and let $(y_i,i\in \mathbb N)$, $y_i< y_{i+1}$, be a homogeneous 
    Poisson point process on $(0,\infty)$ with intensity 
    $p_\varepsilon^M:=\varepsilon^{-\alpha }-M^{-\alpha }$. 
    Consider process $(\mathcal Y_n(u),u\ge 0)$ given by
    \begin{equation}
      \mathcal Y_n(u)=\sum_{i:y_i\le u}\frac{s_n(j)}{t(n)}
    \end{equation}
    or, equivalently, $\mathcal Y_n(u)=\tilde S_n(j)$ for all 
    $u\in [y_{j-1},y_{j})$. We use $\mathcal Y_\varepsilon^M$ to 
    denote the L\'evy 
    process whose  L\'evy measure $\nu_\varepsilon^M$ is given by 
    (a multiple of) the distribution of $s_\infty(j)$,
    \begin{gather}
      \label{e:lm}
      \mathbb E\big[\exp \big(-\lambda \mathcal Y_\varepsilon^M(u)\big)\big]=
      \exp\Big\{ -u\int_0^\infty (1-e^{-\lambda v 
      })\nu_\varepsilon^M(\d v)\Big\},
      \\
      \label{e:nuhh}
      \nu_\varepsilon^M (\d v)= {p_\varepsilon^M} \mathbb 
      P[s_\infty\in \d v]=
      \d v
      \int_\varepsilon^M \alpha z^{-\alpha -2} e^{-v/z}\,\d z.
    \end{gather}
    Condition~\ref{c:conv} implies that as 
    $\mathcal Y_n \to \mathcal K \mathcal Y_\varepsilon^M$ as 
    $n\to \infty$  weakly in the Skorokhod topology  for a.e.~$\boldsymbol \tau $.

    Let $\mathcal R(\mathcal Y_n)$ denote the range of the process 
    $\mathcal Y_n$, 
    $\mathcal R(\mathcal Y_n)=\bigcup_{u\ge 0}\mathcal Y_n(u)$. It 
    follows that
    \begin{equation}
      \mathbb P\big[[a,b]\cap \{\tilde S_n(j), j\in \mathbb N\}=
	\emptyset\big|\boldsymbol \tau \big]=
      \mathbb P\big[ \mathcal R(\mathcal Y_n)\cap [a,b] =\emptyset\big| 
	\boldsymbol \tau \big],
    \end{equation}
    Since the distribution of $s_\infty$ has no atoms, the 
    probability that $a/\mathcal K$ or $b/\mathcal K$ are contained  
    in 
    $\mathcal R(\mathcal Y_\varepsilon^M)$ is zero. Therefore, the 
    weak convergence of $\mathcal Y_n$ implies that 
    \begin{equation}
      \lim_{n\to \infty}
      \mathbb P\big[[a,b]\cap \{\tilde S_n(j), j\in \mathbb N\}
	=\emptyset|\boldsymbol \tau \big]
      =\mathbb P[\mathcal R(\mathcal Y_\varepsilon^M) 
	\cap [a/\mathcal K,b/\mathcal K]=\emptyset | \boldsymbol \tau ].
    \end{equation}
    To estimate the right-hand side of the previous expression
    observe that as $\varepsilon \to 0$ and $M\to \infty$  the L\'evy 
    measure $\nu_\varepsilon^M$ converges to
    \begin{equation}
      \label{e:lmsa}
      \d v \int_0^\infty \alpha z^{-\alpha -2} e^{-v/z}\,\d z= 
      \alpha \Gamma (1+\alpha ) v^{-1-\alpha } \d v.
    \end{equation}
    This is the L\'evy measure of an $\alpha $-stable subordinator. 
    Therefore, as $\varepsilon \to 0$ and $M\to \infty$, the process 
    $\mathcal Y_\varepsilon^M$ converges weakly in the Skorokhod topology to the
    stable subordinator. Using the same reasoning as 
    before we get
    \begin{equation}
      \lim_{\substack{\varepsilon \to 0\\ M\to\infty}}
      \mathbb P[\mathcal R(\mathcal Y_\varepsilon^M) 
	\cap[a/\mathcal K,b/\mathcal K]=\emptyset | \boldsymbol 
	\tau ]= \Asl_\alpha (a/b).
    \end{equation}
    This finishes the proof.
  \end{proof}

  We can now finally estimate the probability of the event 
  $E(n)$ (defined by~\eqref{e:defE}) for large $n$. Fix 
  $\delta >0$.   
  Let $\mathcal B$ be 
  the event that is considered in Condition~\ref{c:large},
  \begin{equation}
    \mathcal B=\Big\{\sum_{i=1}^{\zeta_n}s_n(i)\ge (1+\theta )t(n)\Big\}.
  \end{equation}
  According to Condition~\ref{c:large} we can choose $\xi_n$ such that 
  $\mathbb P[\mathcal B^c|\boldsymbol \tau ]\le \delta$ for 
  a.e.~$\boldsymbol \tau $. 
  We divide the probability space into three disjoint sets:
  \begin{equation}
    \begin{split}
      G_1(n)&=\{\dist(1,\mathcal R(\mathcal Y_n))\le 2 \delta  \text{ or }
	\dist(1+\theta ,\mathcal R(\mathcal Y_n))\le 2 \delta \}\\
      G_2(n)&=\{\dist(1,\mathcal R(\mathcal Y_n))> 2 \delta,  	
	\dist(1+\theta ,\mathcal R(\mathcal Y_n))> 2 \delta  \text{ and } \\
	&\phantom{===}(1,1+\theta )\cap \mathcal R(\mathcal Y_n) \neq \emptyset\}\\
      G_3(n)&=\{ [1-2\delta ,1+\theta +2\delta ]\cap \mathcal R(\mathcal Y_n)=\emptyset\}
    \end{split}
  \end{equation}
  Here, for $A,B\subset \mathbb R$, 
  $\dist(A,B)=\min\{|x-y|:x\in A, y\in B\}$.
  As we have already remarked, this division has the following 
  reasons. Heuristically, on the event $\mathcal B$, to the precision $\delta $ 
  (see Lemma~\ref{l:approx}), any interval that does not intersect 
  $\mathcal R(\mathcal Y_n)$ corresponds to a time period spent in one site of 
  the top and neighbouring shallow traps. On the other hand, the 
  points of $\mathcal R(\mathcal Y_n)$ (or more precisely very short periods 
    preceding them) correspond to times when  no deep trap is 
  visited for a large number of steps. 

  We wish to show that the events  $E(n)$ and $G_3(n)$ are essentially the same. 
  Obviously, for a.e.~$\boldsymbol \tau $ (omitting the 
    conditioning on $\boldsymbol \tau $ in the notation)
  \begin{equation}
    \label{e:suma}
    \mathbb P[E(n)\cap G_3(n)]\le \mathbb P[E(n)]\le
    \mathbb P[G_3(n)]+\mathbb P[G_1(n)]+\mathbb P[E(n)\cap G_2(n)].
  \end{equation}
  We should therefore estimate all quantities in the last display. 
  The probability of $G_1(n)$ is small. Indeed,
  \begin{equation}
    \label{e:metap1}
    \mathbb P[G_1(n)]\le \mathbb P[\dist(1,\mathcal R(\mathcal Y_n))\le2 \delta ]+
    \mathbb P[\dist(1+\theta ,\mathcal R(\mathcal Y_n))\le 2\delta ].
  \end{equation}
  The both probabilities on the right-hand side can be estimated using 
  Lemma~\ref{l:levyconv}, namely it is possible to choose 
  $\varepsilon $ small and $M$, $n$ large enough, such that 
  \begin{multline}
    \mathbb P[\dist(1,\mathcal R(\mathcal Y_n))\le 2\delta ]=
    \mathbb P\big[[1-2\delta ,1+2\delta ]\cap \mathcal R(\mathcal Y_n) = \emptyset \big]\\\le
    \delta + 1 - \Asl_\alpha \Big(\frac{1-2\delta }{1+2\delta }\Big)\le 
      C\delta^{1-\alpha }.
  \end{multline}
  In the same way we estimate the second probability from \eqref{e:metap1}.

  If $\mathcal B\cap G_2(n)$ holds, then there are $j_1<j_2\le\zeta_n$ 
  such that 
  \begin{equation}
    \begin{split}
      \label{e:duk1}
      \tilde S_n(j_1)+2\delta &\le 1 \le \tilde S_n(j_1+1)-2\delta, \\
      \tilde S_n(j_2)+2\delta &\le 1+\theta  \le \tilde S_n(j_2+1)-2\delta,
    \end{split}
  \end{equation}
  and therefore, using Lemma~\ref{l:approx},
  \begin{equation}
    \begin{split}
      \label{e:duk2}
      S_n(r_n(j_1))+\delta t(n) &\le t(n) \le S_n(r_n(j_1+1)) - 
      \delta t(n),\\
      S_n(r_n(j_2))+\delta t(n) &\le (1+\theta )t(n) \le S_n(r_n(j_2+1)) - 
      \delta t(n).
    \end{split}
  \end{equation}
  Hence, according to Condition \ref{c:post}, for all $\varepsilon $ 
  small and $M$ large
  \begin{equation}
    \mathbb P\big[X(t(n))=U_n(j_1),X_n((1+\theta 
	  )t(n))=U_n(j_2)\big|\boldsymbol \tau,\mathcal B\cap G_2 \big] \ge 1-\delta.
  \end{equation}
  Using Condition~\ref{c:noreturn},  it follows  that 
  $\mathbb P[E(n)| G_2(n) \cap \mathcal B]\le \delta$ for all $n$ large 
  enough. Therefore,
  \begin{equation}
    \mathbb P[E(n)\cap G_2(n)]\le \mathbb P[E(n)|G_2(n)\cap 
      \mathcal B]+\mathbb P[\mathcal B^c]\le 2\delta .
  \end{equation}

  At last, we estimate the terms related to $G_3(n)$. 
  Using Lemma~\ref{l:levyconv},
  choosing $\varepsilon $ small and $n$, $M$ large enough, we get
  \begin{equation}
    \big|\mathbb P[G_3(n)]-\Asl_\alpha (1/1+\theta )\big|\le C \delta. 
  \end{equation}
  Finally, by Condition~\ref{c:post} and using a similar reasoning as 
  in \eqref{e:duk1}, \eqref{e:duk2},  
  \begin{equation}
    \mathbb P\big[X(t(n))=X((1+\theta )t(n))\big| \mathcal B\cap G_3\big]\ge 1-\delta .
  \end{equation}
  Since $\delta $ can be taken arbitrarily small, 
  Theorem~\ref{t:main}(i) 
  follows from \eqref{e:suma} and the results of the last three paragraphs.

  The proof of Theorem~\ref{t:main}(ii) proceeds along the same lines  as the
  proof of (i). The only needed changes are the re-definition of $E(n)$ as
  \begin{equation}
    E(n):=\big\{\exists j :
      S_n(r_n(j))\le t(n) < (1+\theta )t(n) \le S_n(r_n(j+1)) \big\}
  \end{equation}
  and the observation that  
  Conditions~\ref{c:post} and~\ref{c:noreturn} are not necessary in 
  this case.
\end{proof}

\begin{proof}[Proof of Proposition~\ref{p:rel}]
  Now we show that  Conditions~(A)--(D) imply 
  Conditions~\ref{c:shallow}--\ref{c:large}. To this end we choose 
  $\xi_n=m r(n)$ with some large constant $m$. The scale $g(n)$ is, 
  of course, the same in  both sets of conditions. 

  (A), (D)$\implies$\textit{Condition~\ref{c:shallow}}. We show that there is 
  a large constant $K$ independent of $\varepsilon $, $M$, and $m$ 
  such that Condition \ref{c:shallow} is satisfied with 
  $h(\varepsilon )=K m \varepsilon^{1-\alpha }$, that is for 
  $\mathbb P$-a.e.~realisation of $\boldsymbol \tau $  and for all 
  large $n$ 
  \begin{equation}
    \label{e:sshallow}
    \mathbb E\Big[\sum_{i=0}^{\xi_n}e_i \tau_{Y_n(i)}\bbone\{Y_n(i)\in 
	T^\varepsilon (n)\}\Big|\boldsymbol \tau \Big]\le 
    K m \varepsilon^{1-\alpha }t(n).
  \end{equation}
  To this end we use the same slicing strategy as in \cite{BCM06,Cer03}.
  We divide $T^\varepsilon(n) $ into disjoint sets 
  $T_{\varepsilon 2^{-i}}^{\varepsilon 2^{-i+1}}(n)$  with 
  $i\in \mathbb N$. We show that there is a large constant $K'$ such that 
  $\mathbb P$-a.s. for all but a finite number of $n$ the following holds:
  for all $i\in \mathbb N$
  \begin{equation}
    \label{e:shahc}
    \mathbb E\Big[
      \sum_{i=0}^{\xi_n}e_i 
      \tau_{Y_n(i)}\bbone\big\{Y_n(i)\in 
	T_{\varepsilon 2^{-i}}^{\varepsilon 2^{-i+1}}
	\big\}\Big|\boldsymbol \tau \Big]\le K'm \varepsilon^{1-\alpha }
    2^{i(\alpha -1)}t(n).
  \end{equation}
  The summation over all $i$ then yields directly the claim \eqref{e:sshallow}.

  Let $p_{n,i}=\mathbb P[x\in T_{\varepsilon 2^{-i}}^{\varepsilon 2^{-i+1}}]$.
  By Condition (A)
  \begin{equation}
    \label{e:pinhc}
    p_{n,i}\le C \varepsilon^{-\alpha }2^{i\alpha }\rho (n)
  \end{equation}
  for some $C$ independent of $n$, $i$, and $\varepsilon $. 
  For all $i\in\mathbb N$ we have
  \begin{equation}
    \begin{split}
      \label{e:sh1}
      \mathbb P\Big[&\mathbb E\Big[
	  \sum_{i=0}^{\xi_n}e_i 
	  \tau_{Y_n(i)}\bbone\big\{Y_n(i)\in 
	    T_{\varepsilon 2^{-i}}^{\varepsilon 2^{-i+1}}
	    \big\}\Big|\boldsymbol \tau \Big]\ge K'm \varepsilon^{1-\alpha }
	2^{i(\alpha -1)}t(n)\Big]\\
      &=
      \mathbb P\Big[\sum_{x\in \mathcal V_n }G^n_{\xi_n}(\boldsymbol 0,x)\tau_x
	\bbone\big\{x\in
	  T_{\varepsilon 2^{-i}}^{\varepsilon 2^{-i+1}}
	  \big\}\ge  K'm \varepsilon^{1-\alpha }
	2^{i(\alpha -1)}t(n)\Big]\\
      &\le 
      \mathbb P\Big[\sum_{x\in \mathcal V_n}G^n_{m r(n)}(\boldsymbol 0,x)
	\bbone\big\{x\in
	  T_{\varepsilon 2^{-i}}^{\varepsilon 2^{-i+1}}
	  \big\}\ge  K'm \varepsilon^{-\alpha }
	2^{i\alpha -1} t(n)g(n)^{-1}\Big].
    \end{split}
  \end{equation}
  By Chebyshev inequality with $\lambda_n$ of Condition (D) this is bounded by
  \begin{equation}
    \label{e:sh1a}
    \le \exp(-\lambda_n m K' \varepsilon^{-\alpha }2^{i\alpha -1} t(n)g(n)^{-1})
    \prod_{x\in \mathcal V_n}\big[1+p_{n,i}\big(e^{\lambda_n G^n_{mr(n)}(\boldsymbol 0,x)}-1\big)\big].
  \end{equation}
  Using $x\ge \log(1+x)$ and \eqref{e:pinhc} we get 
  \begin{equation}
    \le
    \exp\{
      -\lambda_n m K' \varepsilon^{-\alpha }2^{i\alpha -1} t(n)g(n)^{-1}+
      C \varepsilon^{-\alpha } 2^{i\alpha } \rho (n) \sum_{x\in \mathcal V_n}
      \big(e^{\lambda_n G^n_{m r(n)}(\boldsymbol 0,x)}-1\big)\}.
  \end{equation}
  Condition (D) ensures that the last expression is bounded by 
  \begin{equation}
    \le
    \exp\{
      -\lambda_n m \varepsilon^{-\alpha } 2^{i\alpha }t(n)g(n)^{-1}
      ( K' - C  \mathcal K_s) \}
  \end{equation}
  Now it is easy to prove \eqref{e:shahc} and thus  Condition~\ref{c:shallow}: 
  it is sufficient to take $K'> 2 C \mathcal K_s$, to  sum over all $i$ 
  and then to  apply 
  the Borel-Cantelli lemma. This is possible due to the second part 
  of Condition (D).

  \medskip

  (A), (C)$\implies$\textit{Condition~\ref{c:deep}}.
  By the assumptions of Theorem~\ref{t:pot} $\tau_x$ are i.i.d. 
  Therefore, $T_M(n)$ is a Poisson cloud with density that is bounded 
  by $C \rho (n) M^{-\alpha }$. The normalised hitting time of this 
  cloud, $H_n(T_M(n))/r(n)$, is by 
  Condition (C) asymptotically exponentially distributed with mean 
  $ M^{\alpha }/\mathcal K_r C$. So that, by choosing $M$ 
  large, we can make the 
  probability that $H_n(T_M(n))\le \xi_n$ arbitrarily small. 

  \medskip

  (A)--(C)$\implies$\textit{Condition~\ref{c:conv}}. We first prove the 
  following easy lemma which is a consequence of Condition (C) and the
  lack-of-memory property of the exponential distribution
  \begin{lemma}
    \label{l:twoclouds}
    Assume that Condition (C) holds. 
    Let
    $x\in A_n\cup\{\boldsymbol 0\}$, and let $A_n^1$, 
    $A_n^2$ be such that $A_n^1\cup A_n^2=A_n$, 
    $A_n^1\cap A_n^2=\emptyset$ and 
    \begin{equation}
      \lim_{n\to \infty} |A_n^k|/|A_n|=\rho_k/\rho   \qquad k\in\{1,2\}, 
      \rho_k\in(0,1).
    \end{equation}
    Define $H_k:=H_n(A^k_n\setminus\{x\})/r(n)$, $H:=H(A_n\setminus\{x\})/r(n)$. 
    Then, given that $Y_n(0)=x$, the distribution function of 
    $H\bbone\{H_1<H_2\}$ converges uniformly 
    in $x$ to 
    $F(u):=\rho_2/\rho  + \rho_1 (1-e^{- \mathcal K_r \rho  u})/\rho $, $u\ge 0$, 
    and, in particular,
    \begin{equation}
      \label{e:corclaim}
      \lim_{n\to\infty}\sup_{x\in A_n}\Big|
      \mathbb P_x[H_1< 
	H_2]-\frac {\rho_1}\rho  
      \Big | =0.
    \end{equation}
  \end{lemma}

  \begin{proof}
    We define 
    $\nu^k_n(\d u):=\mathbb P_x[H_k\in \d u]$, $k\in\{\emptyset,1,2\}$, and
    \begin{equation}
      f_n(u)=
      \begin{cases}
	\mathbb P_x[H_1<H_2|H=u ] &\text{if 
	  $u\in r(n)^{-1}\mathbb Z$,}\\
	0&\text{otherwise.}
      \end{cases}
    \end{equation}
    Using a decomposition on $H$ and $Y_n(H)$, we 
    get from the Markov property of $Y_n$ for any continuous bounded function $h:\mathbb R^+\to \mathbb R$
    \begin{equation}
      \begin{split}
	\label{e:cora}
	\int h(u) \nu^1_n(\d u) &=
	\int \nu_n(\d u) f_n(u) h(u) \\
	&\quad+ \int \nu_n(\d u) (1-f_n(u))\int \bar \nu_{n,u}(\d v) h(u+v),
      \end{split}
    \end{equation}
    where the measure $\bar\nu_{n,u}$ is defined as
    \begin{equation}
      \bar \nu_{n,u}(\d v)=\sum_{y\in A_n^2\setminus\{x\}}
      \mathbb P_x[Y_n(r(n)H)=y|H=u,H_2<H_1] 
      \mathbb P_y[H_1\in \d v].
    \end{equation}
    Now we take $h(u)=e^{-\eta u}$, $\eta \ge 0$. From 
    Condition (C) it follows that for all $\varepsilon >0$ 
    there is $n_0$ independent of $x\in A_n$ such that for $n>n_0$, $k\in\{\emptyset,1,2\}$,
    \begin{equation}
      \begin{split}
	&\Big | \int \nu_n^k (\d u) e^{-\eta u}-\frac {\mathcal K_r 
	  \rho_k}{\mathcal K_r \rho_k+\eta}\Big| 
	\le \varepsilon, \\
	&\Big | \int \bar\nu_{n,u} (\d u) e^{-\eta u}-
	\frac {\mathcal K_r\rho_1}{\mathcal K_r\rho_1+\eta}\Big| 
	\le \varepsilon. \\
      \end{split}
    \end{equation}
    Inserting this in \eqref{e:cora} and re-arranging it slightly we get
    \begin{equation}
      \label{e:corb}
      \int e^{-\eta u}f_n(u)\nu_n(\d u)=\frac {\mathcal K_r\rho} 
      {\mathcal K_r \rho +\eta }\cdot 
      \frac{\rho_1}{\rho  } + O(\varepsilon ),
    \end{equation}
    from which the first claim of the lemma directly follows. Taking 
    $\eta =0$ in \eqref{e:corb} we get also \eqref{e:corclaim}.
  \end{proof}

  To prove Condition~\ref{c:conv} we need to verify that the sequence 
  $s_n(j)/t(n)$ (see~\eqref{e:defsn} for the definition) 
  converges to the i.i.d.~sequence $s_\infty$.  
  We control first the sequence of depths of 
  visited deep traps:
  \begin{lemma}
    \label{l:deepconv}
    The sequence $(\tau_{U_n(j)}/g(n) ,j\in \mathbb N)$ 
    converges  weakly as $n\to \infty$ to the i.i.d.~sequence 
    $(\sigma_\varepsilon^M(j),j\in\mathbb N)$ defined in \eqref{e:defsigman}. 
  \end{lemma}

  \begin{proof}
    Fix $u\in (\varepsilon ,M)$. Since $\tau_x$ are i.i.d. the set
    $T_\varepsilon^M$ and also 
    its subsets $T_\varepsilon^u$, $T_u^M$ are Poisson clouds with 
    densities converging to 
    $p_\varepsilon^M\rho (n)$, resp. 
    $p_\varepsilon^u\rho (n)$ and $p_u^M\rho (n)$, where 
    $p_a^b:=a^{-\alpha }-b^{-\alpha }$. We can therefore use 
    Condition (C) and  Lemma~\ref{l:twoclouds}:
    uniformly in $x\in T_\varepsilon^M\cup\{\boldsymbol 0\}$ 
    \begin{equation}
      \mathbb P_x\big[
	\tau_{Y_n(H(T_\varepsilon^M\setminus\{x\}))}\le ug(n)\big] 
       =
      \mathbb P_x\big[H(T_\varepsilon^u\setminus\{x\})=H(T_\varepsilon^M\setminus\{x\})\big]
      \xrightarrow{n\to \infty}\frac{ p_\varepsilon^u}{p_\varepsilon^M}.
    \end{equation}
    The lemma then follows easily.
  \end{proof}

  Observe now  that given 
  $U_n(j)=x\in T_\varepsilon^M(n)$, the random variable $s_n(j)/\tau_x$ is 
  exponentially distributed with  mean 
  $G^n_{T_\varepsilon^M\setminus\{x\}}(x,x)=\mathcal K_G f(n)(1+o(1))$ 
  as $n\to\infty$. The error  is uniformly bounded in $x$. 
  Using the Markov property we get
  \begin{equation}
    \begin{split}
      \label{e:tt}
      \mathbb P&[s_n(j)/t(n)\ge u|s_n(1),\dots,s_n(j-1)]\\
      &=\int_{\varepsilon }^M
      \mathbb P[s_n(j)\ge ut(n) | \tau_{U_n(j)}=vg(n), 
	s_n(1),\dots,s_n(j-1)]\\
      &\phantom{=\int_{\varepsilon }^M}
      \times
      \mathbb P[ \tau_{U_n(j)}/g(n)\in \d v|s_n(i),i=1,\dots,j-1]\\
      &=\int_{\varepsilon }^M
      \exp\Big\{-\frac {ut(n)(1+o(1))}{v g(n) \mathcal K_G f(n)}\Big\}
      \mathbb P[ \tau_{U_n(j)}/g(n)\in \d v|s_n(i),i=1,\dots,j-1].
    \end{split}
  \end{equation}
  By definition $g(n)f(n)=t(n)$. The  
  weak convergence of the sequence 
  $\tau_{U_n(j)}/g(n)$ proved in Lemma~\ref{l:deepconv} then yields 
  that the right-hand side of \eqref{e:tt} converges to
  \begin{equation}
    \int_\varepsilon^M e^{-u/\mathcal K_G v} 
    \mathbb P[\sigma_\varepsilon^M(1)\in \d v],
  \end{equation}
  and Condition~\ref{c:conv} is proved.

  \medskip

  (A)--(C)$\implies$\textit{Condition~\ref{c:large}}.
  Fix temporarily $\varepsilon =1/2$, $M=2$.
  Since the distribution of $s_n(i)/t(n)$ converges to the distribution of 
  $s_\infty(i)$, for any $\delta >0$ it is possible to fix a large 
  integer $R$ such that for all $n$ large enough
  \begin{equation}
    \mathbb P\Big[\sum_{i=1}^R s_n(i)\ge (1+\theta )t(n)\Big]\ge 1-\delta /2. 
  \end{equation}
  To satisfy  Condition~\ref{c:large}, the constant $m$ should be fixed such that 
  $\mathbb P[\zeta_n\ge R]\ge 1-\delta /2$. But it can be done easily, 
  since $\zeta_n$ is the number of deep traps that are visited 
  before $\xi_n$, and the distribution of this number converges to 
  the Poisson 
  distribution with mean $p^2_{1/2}m/\mathcal K$ as follows from Condition~(C). 
  Taking now 
  $\varepsilon < 1/2$ or $M>2$, the sum $\sum_{i=1}^{\zeta_n}s_n(i)$ can 
  become only larger. Therefore, for chosen $m$ Condition~\ref{c:large} 
  is verified for all $\varepsilon<1/2$ and $M>2$. 
\end{proof}

\begin{proof}[Proof of Theorem~\ref{t:pot}(i)]
  We want to to verify that Conditions~(A)--(D) imply the weak 
  convergence in the Skorokhod topology  of the process $S( r(n)\cdot)/t(n)$ to an
  $\alpha $-stable subordinator. As usual, it is sufficient to check 
  the convergence of fixed-time distributions and the tightness in 
  $D([0,T])$. 

  Fix $t_0\in (0,T)$. Let $N_n(t_0)$ be the number of deep traps 
  visited in first $t_0 r(n)$ steps,
  \begin{equation}
    N_n(t_0)=\max\{j:r_n(j)\le t_0 r(n)\}.
  \end{equation}
  It follows from Condition~(C) that the distribution of $N_n(t_0)$ 
  converges to the Poisson distribution with mean 
  $t_0 \mathcal K_r p_\varepsilon^M$. Due Proposition~\ref{p:rel} 
  Conditions~\ref{c:shallow}--\ref{c:large} holds, therefore
  using a similar reasoning as to show Lemma~\ref{l:approx}, we can 
  show that for any $\delta $ there are $\varepsilon $ and $M$ such 
  that with probability larger than $1-\delta$ 
  \begin{equation}
    \label{e:afafd}
    t(n)^{-1}\sum_{i=0}^{N_n(t_0)-1} s_n(i) \le t(n)^{-1}S_n( t_0 r(n))\le
    \delta + t(n)^{-1}\sum_{i=0}^{N_n(t_0)}s_n(i).
  \end{equation}
  Now, it is easy to proceed as in Lemma~\ref{l:levyconv} to see that 
  the Laplace transform of 
  $(\mathcal K_G t(n))^{-1}\sum_{i=0}^{N_n(t_0)} s_n(j)$ converges as 
  $n\to \infty$ to \eqref{e:lm} evaluated at $u=\mathcal K_r t_0 $. 
  The distribution with such Laplace transform can be made 
  arbitrarily close (in the weak sense) to the distribution at 
  $\mathcal K_r t_0$ of the 
  $\alpha $-stable subordinator with L\'evy measure \eqref{e:lmsa} by 
  taking $\varepsilon $ small and $M$ large. 
  From this the convergence of the distribution at the fixed time 
  $t_0$ follows. One can get completely analogously the convergence 
  of joint distributions of 
  $(t(n)^{-1}S_n(r_n t_0),\dots,t(n)^{-1}S_n(r_n t_k))$.

  Since $S_n$ is increasing to prove the tightness it is sufficient 
  to check the tightness of the sequence of real random variables 
  $t(n)^{-1}S_n(r(n)T)$. However, this can be done easily using 
  \eqref{e:afafd} with $t_0=T$ and 
  Conditions~\ref{c:shallow}--\ref{c:large}. 
\end{proof}

\section{Aging for short time scales in the Random Energy Model} 
\label{s:rem}

The Random Energy Model (REM) is  the 
simplest mean-field model for  spin-glasses and its static behaviour 
is well understood. The studies of dynamics are much more sparse. The 
first proof of aging in the REM was given in  \cite{BBG03,BBG03b}, 
based on 
renewal theory. The  approach of Section~\ref{s:mechanism}
allows to prove aging on a broad range of shorter time scales. 
We will compare both results later. Before doing 
it, let us define the model and give some motivation why and in what 
ranges of times and temperatures aging occurs. 

The REM is a mean-field model of a spin-glass. It 
consists of $n$ spins that can take values $-1$ or $1$, that is 
configurations of the REM are elements of $\mathcal V_n=\{-1,1\}^n$. The 
energies $\{E_x,x\in \mathcal V_n\}$ of the configurations are 
i.i.d.~random variables. The standard choice of the marginal 
distribution of $E_x$ is centred normal distribution with variance 
$n$. We then define
\begin{equation}
  \tau_x=\exp(\beta  E_x).
\end{equation}

For the dynamics of the REM we require that only one spin can be 
flipped at a given moment. This corresponds to 
\begin{equation}
  \mathcal E_n = \{\langle x,y\rangle\in \mathcal V_n^2:
    \sum_{i=1}^n|x_i-y_i|=2\},
\end{equation}
where  $(x_1,\dots,x_n)$ are the values of individual spins. We use $G_n$ 
to denote the $n$-dimensional hypercube $(\mathcal V_n,\mathcal E_n)$. 
There are many choices for the dynamics of REM, which have the Gibbs 
measure $\boldsymbol \tau $ as a reversible measure. We will 
naturally consider the trap model dynamics \eqref{e:dynamics}. 
which is one of the simplest choices. 
We always suppose that
\begin{equation}
  Y_n(0)=X_n(0)=\boldsymbol 0=(1,\dots,1).  
\end{equation}

We have seen in Section \ref{s:mechanism} that aging occurs only if 
the $\tau_x$ are sufficiently heavy-tailed. However, this certainly 
fails to be true for $\tau_x$ here: an easy calculation gives 
$\mathbb P[\tau_x\ge u]\le u^{-\log u/2\beta^2 n}$ which decreases 
faster than any polynomial. It is therefore clear that if the process 
is given enough time to explore a large part of the configuration 
space and thus to discover the absence of heavy tails, then no aging 
occurs, at least in our picture. On the other hand, on shorter time 
scales the process does not feel the non-existence of heavy tails as 
can be seen from the following simple estimate. Let $\alpha >0$, then
\begin{equation}
  \begin{split}
    e^{\alpha^{2}\beta^2 n/2} 
    \mathbb P\big[\tau_x\ge  u &e^{\alpha \beta^2 n}
      ( \alpha \beta \sqrt{2 \pi n})^{-1/\alpha }
      \big]\\
    &=
    e^{\alpha^{2}\beta^2 n/2} 
    \mathbb P\Big[\frac {E_x}{\sqrt n}\ge \frac
      {\log u + \alpha \beta^2 n -
	\alpha^{-1}\log (\alpha \beta \sqrt{2\pi n})}
      {\beta \sqrt{n}}\Big]\\
    &
    \xrightarrow{n\to\infty} u^{-\alpha },
  \end{split}
\end{equation}
which can be obtained easily using  
$\mathbb P[E_x/\sqrt n\ge u]=(u\sqrt{2 \pi })^{-1}e^{-u^2/2}(1+o(1))$.
Therefore as $n\to \infty$
\begin{equation}
  \label{e:htcomp}
  \mathbb P\Big[\frac {\tau_x}{e^{\alpha \beta^2 n}( \alpha \beta \sqrt{2 \pi n})^{-1/\alpha }}\ge u\Big]=
    e^{-\alpha^2\beta^2n/2}\cdot u^{-\alpha }(1+o(1)).
\end{equation}
In view of the fact that the simple random walk on the hypercube 
almost never backtrack it seems reasonable to presume that 
if the process had time to make only approximately 
$e^{\alpha^2 \beta^2 n/2}$ steps,  then aging could be observed. As 
we will see later this presumption shows to be true.

Before we state the aging result let us remark that there is a much 
stronger relation between ``random exponentials'' $\tau_x$ and 
heavy-tailed random variables. Let $(E_i,i\in \mathbb N)$ be an 
i.i.d.~sequence with the same common distribution as $E_x$.  It 
was proved in \cite{BBM05} that for some properly chosen $Z(n)$ and 
$N(n)$ the normalised sum 
\begin{equation}
  \frac 1 {Z(n)} \sum_{i=1}^{N(n)}e^{\beta \sqrt n E_i}
\end{equation}
converges as $n\to \infty$ in law to an $\alpha $-stable distribution 
with $\alpha $ depending on $\beta$ and  $N(n)$. Our 
methods allow to show that the same is true for the properly 
normalised clock process $S_n$, which is a properly-normalised 
sum of correlated random variables (see \eqref{e:Sdef}).

We now fix  objects for which we verify   
Conditions~\ref{c:shallow}--\ref{c:last}, or more precisely (A)--(D) 
together with Conditions~\ref{c:post},~\ref{c:noreturn}. The scales 
we choose should appear natural in view of~\eqref{e:htcomp}. We define
\begin{align}
  \label{e:hcparmsa}
  t(n)=g(n)&=
  ( \alpha \beta \sqrt{2 \pi n})^{-1/\alpha }\exp(\alpha \beta^2 n),\\
  \label{e:hcparmsr}
  r(n)=\rho (n)^{-1}&=\exp(\alpha^2 \beta^2 n/2),
  \\
  \label{e:hcparmsb}
  \xi_n& = m r(n),  \\
  T_\varepsilon^M(n,\alpha )&=\big\{x\in \mathcal V_n: 
    \tau_x\in\big(\varepsilon g(n) ,Mg(n)\big)\big\}. 
  \label{e:hcparmsc}
\end{align}

\begin{theorem}
  \label{t:hcaging}
  Let the parameters $\alpha \in (0,1)$ and $\beta >0$ be such that 
  \begin{equation}
    \label{e:hcagass}
    3/4 < \alpha^2 \beta^2 / 2\log2 < 1.
  \end{equation}
  Then  for $\mathbb P$-a.e.~random environment $\boldsymbol \tau $
  \begin{equation}
    \lim_{n\to \infty}
    R_n(t(n),(1+\theta )t(n);\boldsymbol \tau ) 
    =\Asl_\alpha (1/1+\theta ).
  \end{equation}
\end{theorem}

\begin{remark}
  1. We believe that the range of the validity \eqref{e:hcagass} of the 
Theorem~\ref{t:hcaging} is not the broadest possible.  The upper 
bound $1$ is correct. If $\alpha^2 \beta^2 / 2\log2>1$, then 
$\xi_n\gg 2^n$.  That means that the state space $\mathcal V_n$ 
becomes too small and the process can feel its finiteness.   
On the other hand, the lower-bound $3/4$ is 
purely technical and can probably be improved. 
It appears because the potential-theoretic methods for the simple 
random walk on the hypercube that we use imply our conditions only 
if the set $T_\varepsilon^M$ is 
sufficiently sparse, or more precisely, if it satisfies the so called minimal 
distance condition (see Proposition~\ref{p:hcpot} and 
  Lemma~\ref{l:poiscloudprop} below).
  
2. Observe that $\alpha $ is a free parameter. It is not fixed by the 
temperature $\beta $. Moreover, the condition \eqref{e:hcagass}  can 
be rewritten as
\begin{equation}
  \alpha^{-1}\beta_c \sqrt{3/4} < \beta  < \alpha^{-1}\beta_c,
\end{equation}
where $\beta_c = \sqrt{2\log 2}$ is the critical temperature in the 
usual REM. This, in particular, means that aging can be observed in 
REM also above the critical temperature, $\beta < \beta_c$. 

\end{remark}

We now compare our results with those of \cite{BBG03b}. To allow this 
comparison we describe very briefly some of the results of this 
paper. In \cite{BBG03b} a discrete-time dynamics is considered, 
however as $n\to \infty$ this dynamics  differs very little from the 
continuous-time dynamics used here. The most important object used 
to prove aging in \cite{BBG03} is a  set of  deep traps
defined by $T_n(E)=\{x\in \mathcal V_n: E_x\ge u_n(E)\}$, 
where the function 
$u_n(E)=\beta_c \sqrt n + E/\beta_c\sqrt n + O(\sqrt{\log n/n})$ is 
chosen in the way that the distribution of $|T_n(E)|$ converges to 
some non-degenerate distribution as $n\to\infty$ and $E$ is kept 
fixed. The mean of this distribution diverges if $E\to -\infty$ 
afterwards. The two-time function considered there is essentially 
the function $R_{T_n(E)}$ averaged over all the starting points in the top, let 
call it $\Pi_n(t_w,t_w+t,E)$. The main aging result of \cite{BBG03b} 
says that for any $\beta > \beta_c$ and  $\varepsilon >0$
\begin{equation}
  \label{e:BBGresult}
  \lim_{t\to\infty}\lim_{E\to -\infty}\lim_{n\to \infty}
  \mathbb P\Big[\Big|\frac{\Pi_n(c_n t, (1+\theta )c_n t, 
	E)}{\Asl_{\beta_c/\beta }(1/1+\theta )} -1 \Big|> \varepsilon \Big]=0,
\end{equation}
where $c_n\sim e^{\beta \sqrt n u_n(E)}$.

Here are the main differences between both results

1. Different two-time functions are considered. We believe that it 
is possible to eliminate the dependence on $T_n(E)$  (i.e.~to convert 
  something  $R_{T_n(E)}$ to something $R_n$) from \eqref{e:BBGresult} 
by a post-processing in the direction of Condition~\ref{c:post}. 
It would be probably more difficult to get the a.s.~convergence instead 
of the convergence in probability. 

2. The main differences are in the considered top sizes and time 
scales. In \cite{BBG03b} the size of the  top is kept bounded as 
$n\to \infty$. This allows to apply ``lumping techniques'' to 
describe the properties of the projection of a simple random walk 
on the hypercube to the top,
that is to prove that that an equivalent of our process $U_n$ (see \eqref{e:Undef}) 
converges to the simple random walk on the complete graph with the 
vertex set $T_n(E)$.  In our case the size of 
the top $T_\varepsilon^M(n)$ increases exponentially with $n$. This 
makes the application of the lumping more difficult. That is why we needed to 
develop different techniques to control the process $U_n$ and in 
particular the random variables $\tau_{U_n(j)}$. These techniques can be found in 
Section~\ref{ss:pthc}. Using them we can verify the crucial 
Conditions~(B) and (C) which imply Condition~\ref{c:conv}. 

3. The time scale 
$c_n\sim e^{\beta \sqrt n u_n(E)}\sim e^{\beta \beta_c n + \beta E/\beta_c}$ 
used in \cite{BBG03b} corresponds to the case 
$\alpha \beta/\beta_c =1$ and is much larger than the scale 
$t(n)\sim e^{\alpha \beta^2n}=e^{\frac{\alpha \beta }{\beta_c}\beta \beta_c n }$. 
These scales become closer if $\alpha \beta /\beta_c$ approaches $1$, 
which is the upper limit of the validity of our theorem. We could 
probably, with some minor notational complications, improve our 
result to $t(n)=e^{\beta \beta_c n} h(n)$ with some $h(n)\to 0$ as 
$n\to \infty$ sufficiently fast, but even in this case $t(n)\ll c_n$. 
Another possibility would be to take the double limit as 
\eqref{e:BBGresult}. This approach may be possibly practicable, 
however,  it does not fall into our framework and we therefore prefer 
not to pursue it.

\subsection{Potential theory on the hypercube}
\label{ss:pthc}
In this section we study properties of the simple random walk on the 
hypercube, in particular hitting times of some relatively diluted but 
still large subsets of the hypercube. 

Let us  introduce some notation first. We write $I(x)$ for the rate 
function of the symmetric Bernoulli distribution on $\{0,1\}$
\begin{equation}
  I(x):=x\log x+(1-x)\log(1-x) +\log 2.
\end{equation}
For $\gamma \in (1/2,1)$ we use $\omega=\omega (\gamma ) $ to denote the unique solution of the equation
\begin{equation}
  \label{e:rovomega}
  I(\omega )=(2\gamma -1)\log 2, \qquad 0\le \omega \le 1/2.
\end{equation}

\begin{proposition}
  \label{p:hcpot}
  Let for all $n\ge 1$ sets $A_n\subset \mathcal V_n$ be such that 
  $|A_n|=\rho_n 2^n$, with ``densities'' $\rho_n$ satisfying 
  $\lim_{n\to \infty}\rho_n2^{\gamma n}=\rho \in (0,\infty)$, for some
  $\gamma \in (1/2,1)$. Let further the sets $A_n$ satisfy the {\em minimal 
    distance condition}
  \begin{equation}
    \label{e:distcond}
    \min\{d(x,y):x,y\in A_n\}\ge (\omega(\gamma )  +\varepsilon )n,
  \end{equation}
  for some small constant $\varepsilon>0$.
  Then for all $s\ge 0$
  \begin{equation}
    \lim_{n\to\infty}\max_{x\in A_n}\Big |
    \mathbb E_x \Big[\exp\Big(-\frac s {2^{\gamma 
	    n}}H\big(A_n\setminus\{x\}\big)\Big)\Big]-\frac \rho {s+\rho }\,\Big|=0.
  \end{equation}
  So that, the hitting time $ H(A_n\setminus\{x\})/2^{\gamma n}$ is asymptotically 
  exponentially distributed  with mean~$1/\rho $.
\end{proposition}

\begin{remark}
  1. Eventually, we will take $A_n$ to be the set $T_\varepsilon^M$, 
  that is a Poisson cloud on $\mathcal V_n$ with  density $\rho(n)$. 
  We will verify that the assumptions of the proposition are 
  a.s.~verified if $\gamma $ is  large enough.  This is the result of 
  Lemma~\ref{l:poiscloudprop} that can be found at the end of this 
  section.
  2. Note that in \cite{BG06} 
  similar results were obtained for sets with the minimal distance 
  between the points $o(n)$, but, on the other hand, the size of these 
  sets should be  much smaller than in our case, namely $O(\log n)$.
\end{remark}

\begin{proof}
  We will use a method introduced by Matthews \cite{Mat88} to show 
  this Proposition. The advantage of this method is that the 
  hitting time of a large set can be controlled by means of much 
  simple estimates on the hitting time of a point.
  
  We define
  \begin{equation}
    \begin{split}
      f^+_n(s)&:=\max \{ \mathbb E_x[\exp (-s 2^{-\gamma n}H(y))]: 
	x,y\in A_n, x\neq y\},\\
      f^-_n(s)&:=\min \{ \mathbb E_x[\exp (-s 2^{-\gamma n}H(y))]: 
	x,y\in A_n, x\neq y\}
  \end{split}
  \end{equation}
  The following lemma is the most important building block of the proof.
  \begin{lemma}[Theorem 1.3 of \cite{Mat88}]
    \label{l:mattthm}
    For any $x\in A_n$
    \begin{equation}
      \label{e:mattthm}
      \begin{split}
	\frac{\Gamma (1/f^-_n(s))} {\Gamma (1/f^+_n(s))}&\cdot
	\frac{\Gamma (|A_n|)} {\Gamma (|A_n|-1)}\cdot
	\frac{\Gamma (|A_n|-2+1/f^+_n(s))} {\Gamma (|A_n|-1+ 1/f^-_n(s))}\\
	\le{}& \mathbb E_x\big[e^{-sH(A_n\setminus \{x\})/2^{\gamma n}}\big]\\
	\le{}&\frac {\Gamma (1/f^+_n(s))} {\Gamma (1/f^-_n(s))}\cdot
	\frac{\Gamma (|A_n|)} {\Gamma (|A_n|-1)}\cdot
	\frac{\Gamma (|A_n|-2+1/f^-_n(s))} {\Gamma (|A_n|-1+ 1/f^+_n(s))},
      \end{split}
    \end{equation}
    where $\Gamma$ is the usual gamma-function.
  \end{lemma}
  \begin{proof}
    The lemma follows directly from Theorem 1.3 of \cite{Mat88}. 
    Using hat for objects as they appear in that paper, we identify
    $\hat N$ with $|A_n|-1$ and $\hat A_i$, $i=1,\dots,\hat N$ 
    with points of $A_n\setminus \{x\}$. Since we are 
    interested in the first visit of the set $A_n\setminus \{x\}$,  the 
    expression \eqref{e:mattthm} is obtained by 
    setting $\hat n=1$ in Matthews' theorem, and by rewriting the 
    products appearing there
    using $\Gamma $ functions.
  \end{proof}
  To apply the previous lemma we need very precise estimates on 
  $f^+_n$ and $f^-_n$. Later in this section we 
  will show 
  \begin{lemma}
    \label{l:fplusminus}
    Let the assumptions of Proposition \ref{p:hcpot} be satisfied.
    Then for all $s> 0$ the functions $f^+_n$ and $f^-_n$ satisfy
    \begin{gather}
      \lim_{n\to \infty}\frac 1 {2^{(1-\gamma )n}f^+_n(s)}=
      \lim_{n\to \infty}\frac 1 {2^{(1-\gamma )n}f^-_n(s)}=s,\\
      \lim_{n\to \infty}\frac 1 {f^-_n(s)}-\frac 1{f^+_n(s)}=0.
    \end{gather}
  \end{lemma}

  To finish the proof of 
  Proposition~\ref{p:hcpot} we will need another elementary technical lemma.
  \begin{lemma}
    \label{l:gammarel}
    Let $a_m$, $b_m$, $\delta_m$ be such that $a_m\to a$, $b_m\to b$ 
    and $\delta_m\to 0$ as $m\to \infty$ with $0<a,b<\infty$. Then
    \begin{equation}
      \label{e:gammarel}
      \lim_{m\to\infty}
      \frac{\Gamma (a_m m)}{\Gamma (a_m m +\delta_m)} \cdot
      \frac{\Gamma (b_m m +\delta_m)}{\Gamma (b_m m) }= 1. 
    \end{equation}
  \end{lemma}
  \begin{proof}
    By Stirling formula, $\Gamma (m)=\sqrt{2\pi }(m-1)^{m-\frac 12}e^{-m+1}(1+o(1))$.
    Therefore, up to a multiplicative correction $1+o(1)$, 
    the expression inside the limit \eqref{e:gammarel} equals
    \begin{equation}
      \frac{\big(1+\frac{\delta_m}{b_m m -1}\big)^{b_m m - \frac 12}}
      {\big(1+\frac{\delta_m}{a_m m -1}\big)^{a_m m - \frac 12}} \cdot
      \frac{(b_m m -1 +\delta_m)^{\delta_m}}
      {(a_m m -1 +\delta_m)^{\delta_m}}.
    \end{equation}
    The logarithm of the numerator of the first fraction satisfies
    \begin{equation}
      \Big(b_m m - \frac 12\Big)\log\Big(1+\frac {\delta_m}{b_m m -1}\Big)\le c 
      \delta_m \xrightarrow{m\to\infty} 0.
    \end{equation}
    The same hold for the denominator. The logarithm of the second 
    fraction is
    \begin{equation}
      \delta_m\log\frac{b_m m -1 +\delta_m}{a_m m -1 +\delta_m}\xrightarrow{m\to\infty} 0.
    \end{equation}
    This finishes the proof.
  \end{proof}
  We now use Lemmas \ref{l:mattthm}--\ref{l:gammarel} to finish the proof 
  of Proposition~\ref{p:hcpot}. Define 
  \begin{equation}
    \begin{aligned}
      m&=m(n)=2^{(1-\gamma )n}\\
      \delta_m&=1/f^-_n(s)-1/f^+_n(s)
    \end{aligned}
    \quad
    \begin{aligned}
      a_m&=2^{-(1-\gamma )n}(|A_n|-2+1/f^+_n(s))\\
      b_m&=2^{-(1-\gamma )n}(1/f^+_n(s))
    \end{aligned}
  \end{equation}
   Then, by 
  Lemma~\ref{l:fplusminus}, $a_m$, $b_m$, $\delta_m$ satisfy the 
  assumptions of Lemma~\ref{l:gammarel}. Therefore, 
  \begin{equation}
    \lim_{n\to\infty}
    \frac {\Gamma (1/f^-_n(s))} {\Gamma (1/f^+_n(s))}\cdot
    \frac {\Gamma (|A_n|-2+1/f^+_n(s))} {\Gamma (|A_n|-2+ 
	1/f^-_n(s))}= 1.
  \end{equation}
  Comparing the last display with the left-hand side of 
  \eqref{e:mattthm} we get using Lemma~\ref{l:fplusminus}
  \begin{equation}
    \label{e:konecmatthm}
    \mathbb E_x\big[e^{-s2^{-\gamma n}H(A_n\setminus\{x\})}\big]\ge
    \frac{|A_n|-1}{|A_n|-2+1/f^-_n(s)}(1+o(1))=\frac \rho 
    {s+\rho }(1+o(1)).
  \end{equation}
  The corresponding upper bound can be obtained analogously. 
\end{proof}
\begin{proof}[Proof of Lemma~\ref{l:fplusminus}]
  To estimate $f^+_n(s)$ and $f^-_n(s)$ we need to compute the 
  Laplace transform 
  $f_n(x,y;s):=\mathbb E_x[\exp(-s2^{-\gamma n}H(y))]$ for all pairs 
  $x,y\in A_n$. This task can be largely simplified using symmetries 
  of the hypercube. First, it is obvious that $f_n(x,y;s)$ depends 
  only on the distance between $x$ and $y$, that is if $d(x,y)=k$, 
  then $f_n(x,y;s)= f_n(k,s)$, where
  \begin{equation}
    f_n(k,s):=\mathbb E_{z_k}[\exp(-s2^{-\gamma n}H(\boldsymbol 0))],
    \qquad z_k:=(\overbrace{-1,\dots,-1}^{k \text{ times}},1,\dots,1).
  \end{equation}
  Second, any walk started at $z_k$ should visit a  point in 
  the distance $l<k$ from $\boldsymbol 0$ before hitting $\boldsymbol 0$. This implies that 
  $f_n(k,s)\le f_n(l,s)$.  Therefore, using the minimal distance condition~\eqref{e:distcond},
  \begin{equation}
    f_n(n,s)\le f^-_n(s)\le f^+_n(s)\le f_n((\omega +\varepsilon)n ,s).
  \end{equation}
  The statement of Lemma~\ref{l:fplusminus} is thus equivalent to 
  \begin{gather}
    \label{e:fN}
    \lim_{n\to \infty}\frac 1 {2^{(1-\gamma )n}f_n(n,s)}=s,\\
    \label{e:rozdil}
    {f_n((\omega +\varepsilon)n ,s)}-{f_n(n, s)}=o(2^{-2(1-\gamma )n}).
  \end{gather}

  We first compute 
  $f_n(n,s)=\mathbb E_{z_n}\big[e^{-s2^{-\gamma n}H(\boldsymbol 0)}\big]$. This 
  computation closely follows Section 3 of \cite{Mat89}. We should be 
  a little bit more careful, since we need to compute the Laplace 
  transform of $H(\boldsymbol 0)$ on a scale that is not typical for it, normally 
  $H(\boldsymbol 0)\sim 2^n$.  (In \cite{Mat89} 
    $\mathbb E_{x}\big[e^{-s2^{-n}H(\boldsymbol 0)}\big]$ was estimated.)   By 
  Fourier methods for random walks on finite groups 
  (\cite{Mat88,Mat89}, see also \cite{Dia88} for the general theory)
  \begin{equation}
    \label{e:furt}
    \mathbb E_x e^{-\lambda H(\boldsymbol 0)}=
    \frac
    {\sum_{y\in \mathcal V_n}(-1)^{x\cdot y}\big[1-e^{-\lambda }
	(1-\frac {2 d(y,\boldsymbol 0)}{n})\big]^{-1}}
    {\sum_{y\in \mathcal V_n}\big[1-e^{-\lambda }
	(1-\frac {2 d(y,\boldsymbol 0)}{n})\big]^{-1}},
  \end{equation}
  where $x\cdot y = \sum_{i=1}^n x_i y_i$ is the standard scalar 
  product in $\mathbb R^n$. Setting $\lambda = s/2^{\gamma n}$, the denominator of 
  \eqref{e:furt} (which does not depend on $x$)  equals 
  \begin{equation}
    \label{e:pt0}
    \sum_{i=0}^n \binom n i \Big[1-e^{-s/2^{\gamma n}}\Big(1-\frac {2i}{n}\Big)\Big]^{-1}.
  \end{equation}
  We expand  $e^{-s/2^{\gamma n}}=1-s/2^{\gamma n}+O(2^{-2\gamma n})$.
  Treating separately the term $i=0$,  the denominator becomes
  \begin{equation}
    \label{e:pt1}
    (1-e^{-s/2^{\gamma n}})^{-1}\bigg\{1+
      \sum_{i=1}^n \frac s {2^{\gamma n}}\cdot \frac n {2i} \binom n i\bigg\} (1+o(1)).
  \end{equation}
  Writing $\frac n{2i}=\sum_{j=0}^{\infty}(1-\frac {2i}n)^j$ for 
  $i\notin\{0,n\}$, we get
  \begin{equation}
    \begin{split}
      \sum_{i=1}^n \frac s {2^{\gamma n}}\cdot \frac n {2i} \binom n i &=
      \frac s {2^{\gamma n}} \sum_{j=0}^\infty \Big(\frac 2 n\Big)^j 
      \sum_{i=1}^{n-1}\binom n i \Big(\frac n2 -i \Big)^j + O(2^{-\gamma n})
      \\&=
      s 2^{(1-\gamma )n}(1+O(n^{-1})).
    \end{split}
  \end{equation}
  To evaluate the inner sum we used the fact that it is, after   
  the normalisation by $2^{-n}$, up to a small error, equal to the $j$-th 
  central moment of the binomial distribution with parameters $n$ and $1/2$. 
  The denominator of \eqref{e:furt} is therefore equal to 
  \begin{equation}
    \label{e:jmenovatel}
    s2^{(1-\gamma )n}(1-e^{-s/2^{\gamma n}})^{-1}(1+o(1))=
    2^n(1+o(1)).
  \end{equation}

  The calculation of the numerator of \eqref{e:furt} for $x=z_n$ can be done 
  analogously. The only difference is that all factors $\binom ni$
  in \eqref{e:pt0} and further
  should be replaced by $C_n(i)=(-1)^i\binom ni$. 
  We need 
  therefore to compute 
  \begin{equation}
    \label{e:numaa}
    \sum_{i=1}^n  C_n(i)\frac s {2^{\gamma n}}\cdot 
    \frac n {2i},
  \end{equation}
  which corresponds to the sum in \eqref{e:pt1}. 

  \begin{lemma}
    \label{l:suma}
    \begin{equation}
      \sum_{i=1}^n (-1)^i \binom ni \frac 1i = -1-\frac 
      12-\dots-\frac 1n.
    \end{equation}
  \end{lemma}
  \begin{proof}
    It is easy to see that
    \begin{equation}
      \frac{\d}{\d x} \big(x^{-1} 
	[(1-x)^n-1]\big)=\sum_{i=1}^n (-1)^i \binom ni \frac {x^i}i.
    \end{equation}
    Therefore,
    \begin{equation}
      \sum_{i=1}^n (-1)^i \binom ni \frac 1i = 
      \int_0^1 \frac {(1-x)^n -1}{x}\,\d x=
      \int_0^1 \frac {v^n -1}{1-v}\,\d v
    \end{equation}
    from which it is easy to finish the proof.
  \end{proof}

  From the last lemma it follows easily that
  \eqref{e:numaa} tends to $0$ as $n\to \infty$, and the numerator of 
  \eqref{e:furt} for $x=z_n$ is equal to 
  \begin{equation}
    \label{e:citatel}
    (1-e^{-s/2^{\gamma n}})^{-1}(1+o(1))=\frac {2^{\gamma n}}{s}(1+o(1)).
  \end{equation}
  Putting together \eqref{e:furt}, \eqref{e:jmenovatel}, and 
  \eqref{e:citatel} we get 
  \begin{equation}
    f_n(n,s)=\mathbb E_{z_n}\big[e^{-s2^{-\gamma n}H(\boldsymbol 0)}\big]=
    \frac 1 {s2^{(1-\gamma )n}}(1+o(1))\qquad \text{as $n\to\infty$}.
  \end{equation}
  This proves \eqref{e:fN}.

  It remains to prove \eqref{e:rozdil}. We use again formula 
  \eqref{e:furt}. Since the denominator does not depend on $x$, we 
  should only compute the numerator of the difference 
  $f_n(h,s)-f_n(n,s)$ (eventually we will take 
    $h=(\omega +\varepsilon )n$ as in \eqref{e:distcond}). 
  Observing that the  terms with $d(\boldsymbol 0,y)=0$ cancel, the 
  numerator of the difference is equal to
  \begin{equation}
    \label{e:cit2}
    \sum_{i=1}^n D_n(i) \Big[
      1-e^{-s/2^{\gamma n}}\Big(1-\frac {2i}n\Big)\Big]^{-1},
  \end{equation}
  where the combinatorial factors $D_n(i)$ are given by
  \begin{equation}
    D_n(i)=\sum_{y:d(\boldsymbol 0,y)=i}\big[(-1)^{z_h\cdot y}-(-1)^i\big].
  \end{equation}
  An easy combinatorial reasoning gives that (using $\tbinom nk:=0$ 
    for $k\notin\{0,\dots,n\}$)
  \begin{equation}
    D_n(i)=-2(-1)^i\sum_{j=1,3,5,\dots}^{n-h} \binom{n-h}j \binom h {i-j}.
  \end{equation}
  Hence, expanding again  $e^{-s/2^{\gamma n}}$, \eqref{e:cit2} can be written as
  \begin{multline}
    -\sum_{i=1}^n 2(-1)^i\sum_{j=1,3,5,\dots}^{n-h}\binom{n-h}j 
    \binom h {i-j}\frac n{2i}\,(1+o(1))\\
    = -n\sum_{j=1,3,5,\dots}^{n-h}\sum_{k=0}^h (-1)^{j+k}
    \binom {n-h} j\binom h k \frac 1{j+k}\,(1+o(1)).
  \end{multline}
  The sum over $k$ can be computed similarly as in 
  Lemma~\ref{l:suma}: for $j\in \mathbb N$
  \begin{equation}
    \sum_{k=0}^h \frac{(-1)^k}{j+k}\binom h k = 
    \int_0^1 x^{j-1}(1-x)^h \,\d x=
    \frac
    {\Gamma(1+h)\Gamma(j)}{\Gamma(1+h+j)}=\frac 1 j {\binom {h+j} j}^{-1}.
  \end{equation}
  Observing that $(-1)^j=-1$, the expression \eqref{e:cit2} equals 
  (up to a multiplicative correction 1+o(1)) 
  \begin{equation}
    n\sum_{j=1,3,5,\dots}^{n-h} \binom {n-h} j
    {\binom {h+j} h}^{-1}\frac 1 j =: n\sum_{j=1,3,5,\dots}^{n-h} K(j).
  \end{equation}
  Taking $h=\lfloor \xi n\rfloor$ for some $\xi\in (0, 1/2)$, and 
  calculating the ratio of the consecutive summands  $K(j)/K(j+2)$, 
  we find that $K(j)$ attains a maximum  
  for $j\sim n(\frac 12 - \xi )$.  The last display  is therefore bounded 
  from above by
  \begin{equation}
    Cn^2 \binom{(1-\xi )n}{(\frac 12 -\xi )n} {\binom{\frac 12 n}{\xi n}}^{-1}.
  \end{equation}
  Standard  Stirling type estimates applied to the previous expression imply that 
  the numerator of the difference $f_n(\xi n,s)-f_n(n,s)$ is, for all $n$ large enough, bounded by 
  \begin{equation}
    \exp[n(I(\xi )+\delta)]=2^{nI(\xi-\delta ' ) /\log 2},
  \end{equation}
  where $\delta, \delta' >0$ can be taken arbitrarily small. Taking now 
  $\xi =\omega +\varepsilon $, choosing $\delta '<\varepsilon $, 
  using  \eqref{e:distcond},  
  the fact that the denominator is of order 
  $2^n$ (see \eqref{e:jmenovatel}), and $I(x)$ is decreasing on $[0,1/2]$,
  it is easy to prove \eqref{e:rozdil}. This completes the proof of 
  Lemma~\ref{l:fplusminus}.
\end{proof}

The next lemma specifies conditions under which the assumptions 
of Proposition~\ref{p:hcpot} are verified for a Poisson cloud on the hypercube. 
It also collects some technical facts about 
this cloud that we will need later.

\begin{lemma}
  \label{l:poiscloudprop}
  Let $A_n$ be a sequence of Poisson clouds on $\mathcal V_n$ with 
  densities $\rho_n$
   satisfying 
  $\lim_{n\to \infty}\rho_n 2^{\gamma n}=\rho \in (0,\infty)$ for 
  some $\gamma \in (0,1)$. Then

  (i) $\mathbb P$-a.s.~ 
  $2^{(\gamma-1)n}|A_n|\in (\rho -\delta ,\rho +\delta )$ 
  for any $\delta >0$ and $n$ large enough.

  (ii) Let $\gamma > 1/2$ and let $\omega '<1/2$ be the unique solution of 
  $ I(\omega ')=2(1-\gamma )\log 2$.
  Then $\mathbb P$-a.s.
  \begin{equation}
    d_\Min:=\min\{d(x,y):x,y\in A_n\}\ge (\omega '-\delta )n
  \end{equation}
  for any $\delta>0 $ and $n$ large enough.

  (iii) If $\gamma > 3/4$, then the set $A_n$ satisfies 
  $\mathbb P$-a.s.~the assumptions of Proposition~\ref{p:hcpot} for 
  all $n$ large enough.   
  
  (iv) The claims (ii), (iii) stay valid if the set 
  $A_n$ is replaced by $A_n\cup \{\boldsymbol 0\}$.
\end{lemma}
As a corollary of Proposition~\ref{p:hcpot} and 
Lemma~\ref{l:poiscloudprop} we get
\begin{corollary}
  \label{c:hcc}
  If $ \alpha^2 \beta^2/2\log 2 >3/4$ then Condition (C) holds for 
  $r(n)$ and $\rho (n)$ as defined in~\eqref{e:hcparmsr} with $\mathcal K_r=1$.
\end{corollary}

\begin{proof}[Proof of Lemma~\ref{l:poiscloudprop}]
  (i) The proof is standard. Using the Chebyshev inequality we get
  \begin{equation}
    \mathbb P\big[\big| 2^{(\gamma -1)n}|A_n|-\rho \big|\ge \delta \big]\le
    \delta^{-2}2^{n(\gamma -1)}.
  \end{equation}
  The Borel-Cantelli lemma then implies the result.

  (ii) We construct the set $A_n$ in the following way. Let 
  $R$ be a binomial random variable with parameters $2^n$ and $\rho_n$, 
  and let $(x_i, i=1,\dots, R)$ be a collection of  randomly chosen 
  points in $\mathcal V_n$, such that given $(x_j, j<i)$ the point $x_i$ 
  is uniformly distributed in $\mathcal V_n\setminus\{x_j:j<i\}$. It 
  is easy to estimate the probability that $x_i$ is too close 
  to some of $x_j$, $j<i$. Indeed, the volume of the ball with radius 
  $(\omega '-\delta )n$ around a point $x$ satisfies
  \begin{equation}
    \big|\{y\in \mathcal V_n:d(x,y)\le (\omega '-\delta )n\}\big|\le 
    e^{n(I(0)-I(\omega '-\delta )+\varepsilon )} 
  \end{equation}
  for all $\varepsilon >0$ and $n$ large enough. Therefore, for all 
  $i\le R$ 
  \begin{equation}
    \mathbb P\big[\min_{j<i}d(x_i,x_j)\le (\omega '-\delta )n\big]\le
    \frac {i-1}{2^n-i+1}e^{n(I(0)-I(\omega '-\delta )+\varepsilon )}.
  \end{equation}
  So that, summing over $i$,
  \begin{multline}
    \mathbb P\big[d_\Min\le(\omega '-\delta )n \big]\le \sum_{r=1}^{2^n}
    \mathbb P[R=r]\sum_{i=1}^r \mathbb P\big[\min_{j<i}d(x_i,x_j)\le (\omega '-\delta )n\big]
    \\ \le \mathbb P\big[R\ge (1+\varepsilon )2^{(1-\gamma )n}\big]
    +(1+\varepsilon )^2 2^{2(1-\gamma )n}e^{n(-I(\omega '-\delta )+\varepsilon )}.
  \end{multline}
  It follows from the assumptions of the lemma that for any $\delta $ 
  we can find $\varepsilon >0$ such that the right-hand side of the 
  last expression decays exponentially with $n$. The 
  Borel-Cantelli lemma then implies (ii).

  (iii) The assumption of Proposition~\ref{p:hcpot} are satisfied if 
  $d_\Min\ge (\omega +\varepsilon )n$ for some $\varepsilon >0$ and 
  $\omega $ defined in \eqref{e:rovomega}. Using part (ii) of this 
  lemma, this condition is satisfied if 
  $\omega '-\delta \ge \omega +\varepsilon $ for some $\varepsilon$, 
  $\delta $. Since $I(\cdot)$ is decreasing on $(0,1/2)$, we get using 
  the definitions of $\omega $, $\omega '$ that $\gamma$ should satisfy 
  $2(1-\gamma )<2\gamma -1$. This is true for any  $\gamma > 3/4$.

  Easy modifications of the proofs of (ii), (iii) to get (iv) are 
  left to the reader.
\end{proof}

\subsection{Proof of Theorem~\ref{t:hcaging}}

We can now finish the proof of 
Theorem~\ref{t:hcaging}. Since Condition (A) is trivially verified 
and (C) follows from Corollary~\ref{c:hcc}, it remains to verify 
Conditions (B), (D), \ref{c:post} and \ref{c:noreturn}  for the choice 
\eqref{e:hcparmsa}--\eqref{e:hcparmsc} of the parameters.

We need one technical lemma first:
\begin{lemma}
  \label{l:blaa1}
  (a) Green's function $G^n_{\xi_n}$ satisfies (recall $\xi_n=mr(n)$)
  \begin{equation}
    \qquad \limsup_{n\to 
      \infty}G^n_{\xi_n}(\boldsymbol 0,\boldsymbol 0)=\limsup_{n\to\infty} G^n_{\xi_n\log \xi_n}(\boldsymbol 0,\boldsymbol 0) =:G_\infty=1.
  \end{equation}
  (b) For all $x\neq \boldsymbol 0$ and for all $n$ large enough 
  \begin{equation}
    G^n_{\xi_n}(\boldsymbol 0,x)\le cn^{-1}.
  \end{equation}
\end{lemma}
\begin{proof}
  (a) Since $G^n_{\xi_n}(x,y)\le G^n_{\xi_n\log \xi_n}(x,y)$ it is 
  sufficient to consider only the second limit.  The probability 
  that $Y_n$ returns to $\boldsymbol 0$ before $\xi_n\log \xi_n$
  tends to $0$ as $n\to \infty$. Indeed,  using the Chebyshev 
  inequality,
  \begin{equation}
    \label{e:blaa1}
    \mathbb P_{z_1}[H(\boldsymbol 0)<\xi_n \log \xi_n]
    \le e^{s \xi_n\log \xi_n/2^n}\mathbb E_{z_1}\big[e^{-sH(\boldsymbol 0)/2^n}\big].
  \end{equation}
  Performing more carefully the same computation as the one on page 138 of 
  \cite{Mat89} it can be proved that   
  \begin{equation}
    \label{e:blaa2}
    \mathbb E_{z_1}[e^{-sH(\boldsymbol 0)/2^n}]=\frac {1+s/n\big[1+O(s2^{-n})+O(1/n)\big]}{1+s[1+O(s2^{-n})][1+1/n+O(n^{-2})]}.
  \end{equation}
  Take now $s=s(n)$ such that $s(n)\to \infty$ and  
  $s(n)\xi_n\log \xi_n/2^n\to 0$ as $n\to \infty$, which is 
  possible by assumption \eqref{e:hcagass} of 
  Theorem~\ref{t:hcaging}. Then  the right-hand side
  of \eqref{e:blaa1}  tends to $0$ as $n \to \infty$.  
  Therefore, for any $\varepsilon >0$ there exists $n_0$ such that for 
  all $n>n_0$   
  $\mathbb P_{z_1}[H(\boldsymbol 0)\le \xi_n\log \xi_n] \le \varepsilon $. 
  Since after every visit of $\boldsymbol 0$ the random walk $Y_n$ 
  should jump to some of its neighbours, we get by iterating the last estimate 
  \begin{equation}
    \mathbb P\Big[\sum_{i=0}^{\xi_n\log \xi_n}\bbone\{Y_n=\boldsymbol 0\}\ge k\Big]\le \varepsilon ^{k-1}.
  \end{equation}
  Therefore 
  \begin{equation}
    G^n_{\xi_n}(\boldsymbol 0,\boldsymbol 0)=
    \mathbb E\Big[\sum_{i=0}^{\xi_n\log \xi_n}\bbone\{Y_n=\boldsymbol 0\}\Big]\le (1-\varepsilon )^{-1}
  \end{equation}
  for all large $n$. Since $\varepsilon $ was arbitrary, the proof 
  of (a) is finished. 

  (b) Since 
  $G^n_{\xi_n}(\boldsymbol 0,x)\le \mathbb P[H(x)\le \xi_n]G^n_{\xi_n}(\boldsymbol 0,\boldsymbol 0)$, 
  and $\mathbb P[H(x)\le \xi_n]$ decreases as 
  $\dist(x,\boldsymbol 0)$ increases, it is sufficient to show that 
  $\mathbb P[H(z_1)\le \xi_n]=\mathbb P_{z_1}[H(\boldsymbol 0)\le \xi_n]\le c n^{-1}$. 
  However, this is a direct consequence  \eqref{e:blaa2} and 
  Chebyshev inequality with $s(n)=n$. 
\end{proof}

\textit{Condition} (D). We verify this condition with  
$\lambda_n=n^{1/2}$. Then, 
\begin{equation}
  \sum_{n=0}^\infty \exp\{-\lambda_n t(n)/g(n)\}=
  \sum_{n=0}^\infty \exp\{-\lambda_n \}< \infty  
\end{equation}
as required. Using the last lemma we get
\begin{equation}
  \sum_{x\in \mathcal V_n}(e^{\lambda_n G^n_{\xi_n}(\boldsymbol 0,x)}-1)=
  e^{2n^{1/2}}+
  \sum_{x\in \mathcal V_n\setminus\{\boldsymbol 0\}}
  (e^{\lambda_n G^n_{\xi_n}(\boldsymbol 0,x)}-1).
\end{equation}
For $x\in \mathcal V_n\setminus \{\boldsymbol 0\}$  
Lemma~\ref{l:blaa1} yields that  $\lambda_n G^n_{\xi_n}(\boldsymbol 0,x)\le c n^{-1/2}$.  
Since $e^u-1\le 2u$ for all $u$ sufficiently close to $0$, the last 
display is bounded from above by 
\begin{equation}
  \label{e:shalow2}
  e^{2n^{1/2}}+
  \sum_{x\in \mathcal V_n\setminus \{\boldsymbol 0\}}
  2G^n_{\xi_n}(\boldsymbol 0,x) \le
  e^{2n^{1/2}}+ 2 \xi (n)\le
  C m r(n) ,
\end{equation}
which was to be proved.

\textit{Condition} (B). This condition follows from the next lemma.
\begin{lemma}
  \label{l:gg}
 Uniformly for $x\in T_\varepsilon^M$ 
  \begin{equation}
    \label{e:hcc1b}
    \lim_{n\to\infty }G^n_{T_\varepsilon^M\setminus\{x\}}(x,x)= 1. 
  \end{equation}
\end{lemma}
\begin{proof}
  Let $H_n'(x)=\min\{i\ge 1: Y_n(i)=x\}$. Then 
  \begin{equation}
    \label{e:hcc2c}
    G^n_{T_\varepsilon^M\setminus\{x\}}(x,x)=\big(\mathbb 
      P_x\big[H_n(T_\varepsilon^M\setminus \{x\})< 
	H_n'(x)\big]\big)^{-1} 
  \end{equation}
  However
  \begin{multline}
    \label{e:hcc2b}
    1-\mathbb P_x\big[H_n(T_\varepsilon^M\setminus \{x\})< H_n'(x)\big]
    \\\le
    \mathbb P_x[H_n'(x)\le \xi_n\log \xi_n] + 
    \mathbb P_x[H_n(T_\varepsilon^M\setminus\{x\})\ge \xi_n \log \xi_n].
  \end{multline}
  The first term is independent of $x$ and converges to $0$, as can 
  be proved e.g.~using Lemma~\ref{l:blaa1}(a). Since 
  $H_n(T_\varepsilon^M\setminus\{x\})/2^{\gamma n}$
  is asymptotically exponentially distributed with mean 
  $1/p_\varepsilon^M$ as follows from Corollary~\ref{c:hcc}, and since $\xi_n=m 2^{\gamma n}$, the second term 
  in \eqref{e:hcc2b} converges also to $0$ uniformly in $x$. The lemma 
  then follows from \eqref{e:hcc2c} and \eqref{e:hcc2b}.
\end{proof}

\textit{Condition~\ref{c:post}. }Let $t_n$ be a deterministic time 
sequence and $A_n$ defined as in \eqref{e:defAn}.  Fix $\varepsilon $ small 
enough such that $h(\varepsilon )$ used in 
Condition~\ref{c:shallow} satisfies $h(\varepsilon )\le \delta^2/4$. 
Then by Chebyshev inequality the time spent in the shallow traps is small,
\begin{equation}
  \mathbb P\Big[\sum_{i=0}^{\xi_n}e_i\tau_{Y_n(i)}\bbone\{Y_n(i)\in 
      T^\varepsilon (n)\}\ge t(n)\delta/2\Big|\boldsymbol \tau \Big]\le \delta/2 .
\end{equation}
Let $B$ be the event from the previous display. Conditionally on $B$ 
and $A_n$, the event 
$\{X_n(t_n)\neq U_n(j_n)\}$  occurs only if $X_n$ returns to $U_n(j_n)$ 
before the time $S_n(r_n(j+1))$. However, the probability of such return 
tends to $0$ by Lemma~\ref{l:blaa1}(a). This finishes the proof.

\textit{Condition~\ref{c:noreturn}. }As follows from (C) the number 
$\zeta_n$ of visited deep traps  converges to Poisson distribution 
with mean $\mathcal K_r m p_\varepsilon^M$. The probability that $Y_n$ 
returns to any of these traps before $\xi_n$ converges to $0$ by 
Lemma~\ref{l:blaa1}(a). From these two facts the condition follows 
easily.

\section{Aging on the torus}
As a second application of the general strategy presented in 
Section~\ref{s:mechanism} we will give a proof of aging on a 
two-dimensional torus. This complements the results of \cite{BCM06} 
about aging on $\mathbb Z^2$. Note that aging on $d$-dimensional 
torus, $d\ge 3$, could be proved analogously, however for  notational 
convenience we treat only the two-dimensional case.

We will consider the following model. Let 
$G_n=(\mathcal V_n, \mathcal E_n)$ be a two-dimensional torus of size 
$2^n$ with nearest-neighbours connection, 
i.e.~$\mathcal V_n=\mathbb Z^2/2^n\mathbb Z^2$, and edge 
$\langle x,y\rangle\in\mathcal E_n$ iff 
$\sum_{i=1}^2 |x_i-y_i| \bmod 2^n = 1$. We use $d(x,y)$ to denote the 
graph distance of $x,y\in \mathcal V_n$. Let further 
$\boldsymbol \tau =\{\tau^n_x\}$, $x\in \mathcal V_n$, 
$n\in \mathbb N$, be a collection of positive i.i.d.~random variables 
satisfying 
\begin{equation}
  \label{e:tauZ}
  \mathbb P[\tau_x^n\ge u]=u^{-\alpha }(1+o(1))\qquad ( u\to\infty). 
\end{equation}
For simplicity we assume that $\mathbb P[\tau_x\ge 1]=1$. 
Given graph $G_n$ and the random environment 
$\boldsymbol \tau $ we will consider the Markov processes $X_n$ defined 
in Section~\ref{s:introduction}. For this process we will show:

\begin{theorem}
  \label{t:Zdaging}
  Let $t(n)=2^{2n/\alpha }n^{1-(\gamma /\alpha)}$ for some
  $\gamma \in (0,1/6)$. Then 
  for $\mathbb P$-a.e.~realisation of the random environment 
  $\boldsymbol \tau $
  \begin{equation}
    \lim_{n\to\infty} R_n(t(n),(1+\theta )t(n);\boldsymbol \tau )
    =\Asl_\alpha (1/1+\theta ).
  \end{equation}
\end{theorem}

To show Theorem \ref{t:Zdaging} we will verify that  
Conditions~(A)--(D), \ref{c:post}, and \ref{c:last} hold in our case for the 
parameters
\begin{equation}
  \label{e:parmZ}
  \begin{gathered}
    t(n)=2^{2n/\alpha }n^{1-\gamma /\alpha },\qquad\qquad
    \xi_n =m r(n)= m 2^{2n}n^{1-\gamma },\\
    g(n)=\rho (n)^{-1/\alpha }=2^{2n/\alpha }n^{-\gamma /\alpha },\qquad\qquad
    f(n)=n,\\
    T_\varepsilon^M(n)=
    \{x\in \mathcal V_n:
      \tau_x\in (\varepsilon g(n) ,M g(n))\}.
  \end{gathered}
\end{equation}

The main motivation for Theorem~\ref{t:Zdaging} is to extend the 
range of aging scales on $\mathbb Z^2$ and mainly to really explore 
the extreme values of the random landscape. Namely, the BTM on the whole 
lattice $\mathbb Z^2$ does not find the deepest traps that are close to its 
starting position. In the first $2^{2n}$ steps, it gets to the distance 
$2^{n}$ and visits $O(2^{2n}/\log(2^{2n}))=O(2^{2n}/n)$ sites. 
Therefore, the deepest visited trap has a depth of order 
$2^{2n/\alpha }/n^{1/\alpha }$, which is much smaller that the depth 
of the deepest trap in the disk with radius $2^n$, that is 
$2^{2n/\alpha }$. Eventually, the process visits also this deepest 
trap, however it will be too late. This trap will no longer be 
relevant for the time change since much deeper traps will have 
already be 
visited. The deepest trap is relevant only if the random walk stays 
in the neighbourhood of its starting point a long enough time. One 
way to force it to stay is to change $\mathbb Z^2$ to the torus. By 
changing the size of the torus relatively to the number of considered 
steps, i.e.~by changing $\gamma $, different depth scales become 
relevant for aging.  

The range of possible values of $\gamma \in (0,1/6)$ has, as in the 
REM case, a natural bound and an artificial one. It is natural that 
$\gamma <0$ cannot be considered, since in this case the number of 
steps $\xi_n$ of  the simple walk is larger than 
$2^{2n}(\log2^n)^{1+\varepsilon }$, $\varepsilon >0$. 
Therefore, its occupation probabilities are very 
close to the uniform measure on the torus, that is the process is 
almost in equilibrium. The other bound, $\gamma =1/6$, comes from the 
techniques that we use. We do not believe it to be meaningful since we 
expect the theorem to hold for all $\gamma >0$. Actually, the result 
for  $\gamma >1$ follows easily from the result of \cite{BCM06} for 
the whole lattice.  In this case the size of the torus is much larger 
than $\xi_n^2$. So that, the process has no time  to discover 
the difference between the torus and $\mathbb Z^2$. We also know that 
Theorem~\ref{t:Zdaging} holds also in the window $[1/6,1]$ since 
it can be proved by the same methods as for  $\mathbb Z^2$, 
\cite{BCM06}.  However, to keep the presentation in this paper 
compact and to avoid the more complicated coarse-graining 
techniques of \cite{BCM06}, we prefer to stick to 
$\gamma \in (0,1/6)$.

To show our conditions we need again 
a little bit of the potential theory for the simple random walk $Y_n$ on 
the torus $G_n$.

\subsection{Potential theory on the torus}
The following proposition is an equivalent of Proposition~\ref{p:hcpot}. 

\begin{proposition}
  \label{p:Zdpot}
  Let $A_n\subset \mathcal V_n$ be such that $|A_n|=\rho_n 2^{2n}$ 
  with the ``density'' $\rho_n$ satisfying 
  $\lim_{n\to\infty} 2^{2n}n^{-\gamma }\rho_n=\rho$ 
  for some $\gamma \in (0,1/6)$ and $\rho \in (0,\infty)$. 
  Let further $A_n$ satisfy the minimal distance condition
  \begin{equation}
    \min\{d(x,y):x,y\in A_n\}\ge 2^{n}n^{-\kappa }=:r_\Min,
  \end{equation}
  for some $\kappa >0 $. Then,
  for $\mathcal K_r=\pi (2\log 2)^{-1}$,
  \begin{equation}
    \lim_{n\to\infty}\max_{x\in A_n}
    \bigg|\mathbb E_x
    \Big[\exp\Big(-\frac s {2^{2n}n^{1-\gamma }}
	H(A_n\setminus\{x\})\Big) \Big]- 
    \frac {\mathcal K_r\rho }{s+\mathcal K_r\rho }\bigg|=0.
  \end{equation}
\end{proposition}
\begin{proof}
  To proof this proposition we use again the methods from 
  \cite{Mat88}, namely the formula~\eqref{e:mattthm}. To apply it, we 
  need to get precise estimates of the following functions (we use 
    translation invariance of the torus and the minimal distance 
    condition to simplify the definitions)
  \begin{equation}
    \begin{split}
      f^+_n(s)&=\sup_{x\in \mathcal V_n^\out}
      \mathbb E_x\Big[\exp\Big(
	  -\frac s {2^{2n}n^{1-\gamma }}H(\boldsymbol 0)\Big) \Big]\\
      f^-_n(s)&=\inf_{x\in \mathcal V_n^\out}
      \mathbb E_x\Big[\exp\Big(
	  -\frac s {2^{2n}n^{1-\gamma }}H(\boldsymbol 0)\Big) \Big],
    \end{split}
  \end{equation}
  where 
  $\mathcal V_n^\out=\{x\in \mathcal V_n :d(\boldsymbol 0,x)\ge r_\Min\}$, 
  and $\boldsymbol 0=(0,0)$. 
  Using the same reasoning as between \eqref{e:mattthm} and 
  \eqref{e:konecmatthm} it is easy to show that 
  Proposition~\ref{p:Zdpot} follows from
  \begin{gather}
    \label{e:tcl1}
    \lim_{n\to\infty}\frac 1 {n^{\gamma }f^{\pm}_n(s)}=\mathcal K_r^{-1}s,\\
    \label{e:tcl2}
    \lim_{n\to\infty}\Big(\frac 1 {f^-_n(s)}-\frac 1 {f^+_n(s)}\Big)=0.
  \end{gather}

  To estimate the functions $f^+_n$ and $f^-_n$ we use methods 
  developed in \cite{Cox89,CD02}. Following these papers we 
  denote by $q_k(x,y)$ for the (smoothed) probability that a simple random walk 
  on $\mathbb Z^2$ started at $x$ is at $y$ after $k$ steps, more 
  precisely we define
  \begin{equation}
    \label{e:smooth}
    q_k(x,y):=\big\{ \mathbb P_x[Y(k)=y]+ \mathbb P_x[Y(k+1)=y]\big\}/2
    \qquad x,y\in \mathbb Z^2,
  \end{equation}  
  where $Y$ is a simple random walk on $\mathbb Z^2$.
  Let $L:=2^n$. We use
  \begin{equation}
    q_k^n(x,y):=\sum_{z\in \mathbb Z^2}q_k(x,y+zL),  \qquad x,y\in \mathcal V_n,
  \end{equation}
  to denote the transition probability of the simple random walk on the torus 
  $\mathcal V_n$. Let 
  \begin{equation}
    \label{e:defGreen}
    G^n(x;s)=\sum_{k=0}^{\infty}e^{-sk/h(n)} q^n_k(\boldsymbol 0,x)
  \end{equation}
  with $h=h(n)=2^{2n}n^{1-\gamma }$. 
  As is shown e.g.~in \cite{Cox89}, p.~1339, 
  \begin{equation}
    f^\pm_n(s)=\sup_{x\in \mathcal V_n^\out} (\inf_{x\in \mathcal 
	V_n^\out}) \frac{G^n(x;s)}{G^n(\boldsymbol 0;s)}.
  \end{equation}
  Adapting slightly the calculation in \cite{Cox89}, p.~1340, we find 
  that the common denominator in the previous formula satisfies as 
  $n\to \infty$
  \begin{equation}
    \label{e:tcl2a}
    G^n(\boldsymbol 0;s)=\frac {2n} \pi \log 2 (1+o(1)).
  \end{equation}
  In view of this, \eqref{e:tcl1}, \eqref{e:tcl2} are equivalent to 
  \begin{gather}
    \label{e:tcl3}
    \sup_{x\in \mathcal V_n^\out} (\inf_{x\in \mathcal V_n^\out}) 
    G^n(x;s)=\frac {n^{1-\gamma }}s \,  (1+o(1)),\\
    \label{e:tcl4}
    \sup_{x\in \mathcal V_n^\out} 
    G^n(x;s)-\inf_{x\in \mathcal V_n^\out}G^n(x;s)=o(n^{1-2\gamma }).
  \end{gather}
  To verify these two claims we need very precise estimates on 
  $q_k^n(\boldsymbol 0,x)$. Fortunately, the estimates from 
  \cite{CD02} could be used with a small refinement: 

  \begin{lemma} [\cite{CD02} Lemma 3.1]
    \label{l:CD}
    (a) Let $\varepsilon_n=1/\sqrt{\log L}$. There is a finite 
    constant $C$ such that 
    \begin{equation}
      \sup_{k\ge \varepsilon_n L^2}
      \sup_{x\in \mathcal V_n} \varepsilon_n L^2 \rho_k^n(\boldsymbol 0,x)\le C.
    \end{equation}

    (b) If $s_n\to \infty$ as $n\to \infty$, then
    \begin{equation}
      \sup_{k\ge s_n L^2}
      \sup_{x\in \mathcal V_n}L^2|\rho_k^n(\boldsymbol 0,x)-L^{-2}|=o(s_n^{-1}).
    \end{equation}

    (c) If $u_n\to \infty$ as $L\to\infty$, then
    \begin{equation}
      \sup_{x\in \mathcal V_n}\sup_{u_n(1+|x|)^2\le k \le 
	\varepsilon_n L^2}
      |\pi  k \rho_k^n(\boldsymbol 0,x)-1|\to 0.
    \end{equation}

    (d) There is a finite constant $C$ such that
    \begin{equation}
      \sup_{k\ge 0,\,x\in \mathcal V_n}(1+|x|^2)\rho_k^n(\boldsymbol 0,x)\le C.
    \end{equation}
  \end{lemma}
  \begin{proof}
    The claims (a), (c) and (d) are proved in \cite{CD02}. The claim 
    (b) is a refinement of \cite{CD02}, where the same 
    expression was proved to be $o(1)$.   As in \cite{Cox89} we use a very 
    precise expansion of $q_k(\boldsymbol 0,x)$ from \cite{BR-R76}, Corollary 22.3, 
    to prove  (b). As stated there 
    \begin{equation}
      q_k(\boldsymbol 0,x)=\phi_k(x)+\psi_k(x).
    \end{equation}
    Here $\phi_k(x)=(\pi k)^{-1}\exp(-|x|^2/k)$ with 
    \begin{equation}
      |x|=\min\{\sqrt{(x_1+z_1L)^2+(x_2+z_2L)^2}:z\in \mathbb Z^2\},
    \end{equation} 
    and 
    \begin{equation}
      \label{e:defpsi}
      \psi_k(x)=\phi_n(x)\sum_{r=1}^w k^{-r/2}B_r(x/\sqrt k) + e_w(x,k),
    \end{equation}
    where $w$ is an integer larger than $2$, each $B_r(x)$ is a polynomial of degree at most $3r$, and 
    \begin{equation}
      \label{e:erbound}
      \sum_{x\in \mathbb Z^2}|e_w(x,k)|=o(k^{-w/2}) \qquad(n\to\infty).
    \end{equation}
    To get (b), it is sufficient to improve estimates (i), (ii) on 
    page 1373 of \cite{CD02}, namely  to show
    \begin{gather}
      \label{e:tcl5}
      \textstyle
	L^2 \sum_{z\in \mathbb Z^2}\phi_{[s_nL^2]}(x+Lz)-L^2=o(s_n^{-1}),\\
      \textstyle
      \label{e:tcl6}
	L^2 \sum_{z\in \mathbb Z^2}\psi_{[s_nL^2]}(x+Lz)=o(s_n^{-1}).
    \end{gather}

    Let $x'=x/L$, then  the left-hand side of \eqref{e:tcl5} can be 
    written as 
    \begin{multline}
      \label{e:huhu}
      \sum_{z\in \mathbb Z^2}\frac 1{s_n\pi }
      \exp\Big\{-\frac{|x'+z|^2}{s_n}\Big\}
      -\int_{\mathbb R^2}\frac1{s_n \pi }
      \exp\Big\{-\frac{|x'+y|^2}{s_n}\Big\}\, dy
      \\=
      \frac 1 \pi 
      \sum_{z\in \mathbb Z^2/\sqrt{s_n}}
      \int_{I_{s_n^{-1/2}}}
      e^{-|x''+z|^2}
      - e^{-|x''+z+y|^2}\, dy,
    \end{multline}
    where $x''=x'/\sqrt{s_n}$ and 
    $I_{\delta }=[-\delta /2,\delta /2]^2$. Let $h(u)=e^{-|u|^2}$. 
    It is easy to show using Taylor expansion up to forth order and 
    eliminating the odd terms  that 
    \begin{equation}
      \begin{split}
	\begin{split}
	  \int_{I_\delta }h(u+v)-h(v)\,\d v &= 
	  \int_{I_\delta }\Big[
	  \frac 1 2 
	  \Big(\partial_{11} h(u) v_1^2 + \partial_{22}h(u)v_2^2\Big) +
	  R(u)\Big(\frac{\delta}{2}\Big)^4\Big] \d v,
	  \\&= 
	  \frac {\delta^4}8(\partial_{11} h(u) + 
	    \partial_{22}h(u))+R(u)\delta^6.
	    \d v
	\end{split}
      \end{split}
    \end{equation}
    where the reminder  $R(u)$ is bounded by the sum of suprema 
    of all derivatives of forth order over  the set $u+I_\delta $. 
    Since 
    \begin{equation}
      \sum_{z\in \mathbb Z^2/\sqrt{s_n}}\partial_{11} h(z)+\partial_{22} h(z) +R(z)<\infty
    \end{equation}
    we get
    that \eqref{e:huhu} is bounded from above by 
    $C(\sqrt{s_n})^{-4}=Cs_n^{-2}$ uniformly in $x\in \mathcal V_n$. 
    This proves~\eqref{e:tcl5}.

    To show \eqref{e:tcl6} is much simpler. It follows from the
    calculation at page 1343 of \cite{Cox89} that the contribution 
    of the first part of $\psi $ (see \eqref{e:defpsi}) is actually 
    bounded by $C L^{-1}\ll s_n^{-1}$. Similarly, \eqref{e:erbound} 
    yields that the contribution 
    of the second part of $\psi $  is bounded by
    \begin{equation}
      L^2\sum_{z\in \mathbb Z^2}|e_w(x+zL,s_nL^2)|\le
      L^2\sum_{z\in \mathbb Z^2}|e_w(z,s_n L^2)|=L^2o(L^{-2}s_n^{-w/2}).
    \end{equation}
    It is now sufficient to take $w=3$ to finish the proof of 
    \eqref{e:tcl6} and therefore of Lemma~\ref{l:CD}.
  \end{proof}

  We can now finish the proof of \eqref{e:tcl3} and \eqref{e:tcl4}. 
  Choose $s_n$ such that 
  $n^\gamma \ll s_n \ll n^{1/2-2\gamma} \ll n^{1-\gamma }$, which is 
  possible for all $\gamma \in (0,1/6)$. Then, 
  first, by Lemma~\ref{l:CD}(b), uniformly in $x$,
  \begin{equation}
    \begin{split}
      \label{e:oo1}
      \sum_{k=s_nL^2}^\infty &e^{-sk/h}\rho_k^L(\boldsymbol 0,x)=
      \sum_{k=s_nL^2}^\infty e^{-sk/h}L^{-2}(1+o(s_n^{-1}))\\
      &=
      L^{-2}e^{-ss_n/n^{1-\gamma}}
      (1-e^{-s/h})^{-1}(1+o(s_n^{-1}))\\
      &=
      \frac {n^{1-\gamma }}s(1+o(n^{-\gamma })).
    \end{split}
  \end{equation}
  Second, by Lemma~\ref{l:CD}(a), uniformly in 
  $x\in \mathcal V^\out_n$ (recall $h=2^{2n}n^{1-\gamma }$),
  \begin{equation}
    \label{e:oo2}
    \sum_{k=\varepsilon_nL^2}^{s_nL^2} e^{-sk/h}\rho_k^L(\boldsymbol 0,x)\le
    C e^{-sk/h} \frac C{\varepsilon_L L^2}\le \frac 
    {Cs_n}{\varepsilon_n s}=o(n^{1-2\gamma }).
  \end{equation}
  Third, by Lemma~\ref{l:CD}(c), 
  \begin{equation}
    \label{e:oo3}
    \sum_{k=u_n(1+|x|)^2}^{\varepsilon_nL^2}\!\!\!\! e^{-sk/h}\rho_k^L(\boldsymbol 0,x)\le
     \!\!\!\sum_{k=u_n(1+|x|)^2}^{\varepsilon_nL^2} \frac Ck\le
    C \Big[\log \frac{\varepsilon_n} {u_n} -2\log\frac{1+|x|}{L}\Big].
  \end{equation}
  Since $1\ge(1+|x|)/L>n^{-\kappa }$ for $x\in \mathcal V^\out_n$,  
  choosing $u_n=n^{1-3\gamma}$, the last expression is bounded by 
  $C'\log n\ll n^{1-2\gamma }$. Finally,
  \begin{equation}
    \label{e:oo4}
    \sum_{k=0}^{u_n(1+|x|)^2}e^{-sk/h}\rho_k^L(\boldsymbol 0,x)\le
    \sum_{k=0}^{u_n(1+|x|)^2}C (1+|x|^2)^{-1}\le C u_n \ll n^{1-2\gamma }.
  \end{equation}
  Putting together \eqref{e:oo1}--\eqref{e:oo4} we get that for all 
  $x\in \mathcal V_n^\out$
  \begin{equation}
    G^n(x;s)=\frac{n^{1-\gamma }}s(1+\bar o(1))+o(n^{1-2\gamma }),
  \end{equation}
  which is equivalent to \eqref{e:tcl3}, \eqref{e:tcl4}. (Here again 
    $\bar o$ denotes an error that is independent of 
    $x\in \mathcal V_n^\out$.) This finishes the proof of 
  Proposition~\ref{p:Zdpot}.
\end{proof}

As in the case of the REM,  the assumptions
of Proposition~\ref{p:Zdpot} are verified for a Poisson cloud on the torus: 
\begin{lemma}
  \label{l:poiscloudpropZd}
  Let $A_n$ be  
  Poisson clouds  with densities $\rho_n$
  satisfying 
  $\lim_{n\to \infty}\rho_n 2^{2n}n^{-\gamma }=\rho \in (0,\infty)$ for 
  some $\gamma \in (0,1)$. Then

  (i) $\mathbb P$-a.s.~ 
  $n^{-\gamma }|A_n|\in (\rho -\delta ,\rho +\delta )$ 
  for any $\delta >0$ and $n$ large enough.

  (ii) Let $\gamma \in (0,1)$. 
  Then for any $\kappa>1+\gamma $, $\mathbb P$-a.s. for $n$ large enough
  \begin{equation}
    d_\Min:=\min\{d(x,y):x,y\in A_n\}\ge 2^nn^{-\kappa  }.
  \end{equation}

  (iii) If $\gamma \in (0,1/6)$, then the set $A_n$ satisfies 
  $\mathbb P$-a.s.~the assumptions of Proposition~\ref{p:Zdpot} for 
  all $n$ large enough.    
  
  (iv) The claims (ii), (iii) stay valid if the set 
  $A_n$ is replaced by $A_n\cup \{\boldsymbol 0\}$.
\end{lemma}
\begin{proof}
  The proof follows the same lines as the proof of 
  Proposition~\ref{p:hcpot} and is left to the reader.
\end{proof}
\begin{corollary}
  If $\gamma \in (0,1/6)$ then Condition (C) holds for $r(n)$, 
  $\rho (n)$ from \eqref{e:parmZ}.
\end{corollary}

\subsection{Proof of Theorem~\ref{t:Zdaging}}
\label{ss:c1}
It remains to verify Conditions (B), (D), \ref{c:post} and 
\ref{c:noreturn} for objects in \eqref{e:parmZ}.  The condition (A) 
follows trivially from \eqref{e:tauZ}.

\textit{Condition} (D).
 We need first an  estimate on the Green's function 
$G^n_{\xi_n}(0,x)$. Since (up to a very small error due to 
  smoothing \eqref{e:smooth} that we can ignore) 
$G^n_{\xi_n}(0,x)=\sum_{k=0}^{\xi_n}q^n_k(0,x)$, we use 
Lemma~\ref{l:CD}. First, by (b) of this lemma,
\begin{equation}
  \label{e:starj}
  \sum_{k=s_nL^2}^{\xi_n}q^n_k(0,x)\le
  \sum_{k=s_nL^2}^{\xi_n}CL^{-2}\le C\xi_nL^{-2}\le Cmn^{1-\gamma }.
\end{equation}
Second, by (a),
\begin{equation}
  \sum_{k=\varepsilon_nL^2}^{s_n L^2}q^n_k(0,x)\le 
  (s_n-\varepsilon_n)L^2\frac{C}{\varepsilon_nL^2}\le Cs_n 
  \varepsilon_n^{-1}\ll n^{1-\gamma }
\end{equation}
if $s_n$ is such that $s_n\to\infty$ and $s_n\ll n^{1/2-\gamma }$. 
Further, by (c),
\begin{equation}
  \sum_{k=u_n(1+|x|)^2}^{\varepsilon_nL^2}q^n_k(0,x)\le 
  \sum_{k=u_n(1+|x|)^2}^{\varepsilon_nL^2}\frac Ck\le
  C\log\frac{L^2\varepsilon_n}{u_n(1+|x|)^2},
\end{equation}
and by (d)
\begin{equation}
  \label{e:stard}
  \sum_{k=0}^{u_n(1+|x|)^2}q^n_k(0,x)\le u_n. 
\end{equation}
Taking $u_n=n^{1-\gamma }$ and putting together 
\eqref{e:starj}--\eqref{e:stard} we get that for all 
$x\in \mathcal V_n$
\begin{equation}
  G_{\xi_n}^n(0,x)\le C_1 m n^{1-\gamma }+C_2 \log\frac{L^2\varepsilon_n}{u_n(1+|x|)^2}.
\end{equation}

We can now prove Condition (D) for $\lambda_n=n^{\delta -1}$ with 
$0<\delta < \gamma $. We first treat $x\in \mathcal V_n$ 
with $d(0,x)\ge 2^n n^{-\kappa }$, $\kappa >1$. In this case 
$G_{\xi_n}^n(0,x)\le C_1 m n^{1-\gamma }$ and thus 
$\lambda_nG_{\xi_n}^n(0,x)\ll 1$. Therefore
\begin{equation}
  \sum_{x\in \mathcal V_n:|x|\ge 2^nn^{-\kappa }}
  (e^{\lambda_n G^n_{\xi (n)}(0,x)}-1)\le
  \sum_{x\in \mathcal V_n:|x|\ge 2^nn^{-\kappa }}
  2\lambda_n G^n_{\xi (n)}(0,x)\le 2 \lambda_n m r(n).
\end{equation}
The rest of the sum is negligible. Indeed, 
\begin{equation}
  \begin{split}
    \sum_{x\in \mathcal V_n:|x|\le 2^nn^{-\kappa }}&
    e^{\lambda_n G^n_{\xi (n)}(0,x)}
    \le
    e^{C_1 m n^{\delta -\gamma }}
    \sum_{|x|\le 2^nn^{-\kappa }}
    \Big(\frac{c 2^{2n}n^{1/2-\gamma }}{(1+|x|^2)}\Big)^{C_2 \lambda_n}
    \\
    &\le
    (c 2^{2n}n^{1/2-\gamma })^{C_2\lambda_n}
    \int_{1}^{2^nn^{-\kappa }}
    z^{1-2 C_2 \lambda_n}\, \d z\\
    \\&\le 
    C(c 2^{2n}n^{1/2-\gamma })^{C_2\lambda_n}
    (2^{n}n^{-\kappa })^{2-2 C_2 \lambda_n}\le C 2^{2n}n^{-2\kappa 
    }\ll \lambda_n r(n).
  \end{split}
\end{equation}

\textit{Condition} (B). This condition follows from the next lemma.
\begin{lemma}
  \label{l:ggZd}
  Uniformly for $x\in T_\varepsilon^M$ 
  \begin{equation}
    \label{e:hcc1bZd}
    \lim_{n\to\infty }n^{-1}G^n_{T_\varepsilon^M\setminus\{x\}}(x,x)= 
    \frac {2\log2} \pi . 
  \end{equation}
\end{lemma}
\begin{proof}
  By Lemma~\ref{l:poiscloudpropZd} there is no point of the top in the 
  disk $D_x(2^nn^{-\kappa })$ with radius $2^nn^{-\kappa }$ around 
  $x\in T_\varepsilon^M$. Therefore, by e.g.~Theorem~1.6.6 of \cite{Law91},
  \begin{equation}
    G^n_{T_\varepsilon^M\setminus\{x\}}(x,x)\ge
    G^n_{D_x(2^nn^{-\kappa })^c}(x,x)=
    \frac 2 \pi \log (2^nn^{-\kappa 
    })+ O(1) = \frac{ 2\log 2} \pi n+o(n). 
  \end{equation}
  Further, after hitting the boundary of $D_x(2^{n}n^{-\kappa })$ the 
  simple random walk $Y_n$ has probability of order $n^{-\gamma }$ to 
  return back to $x$ before hitting $T_\varepsilon^M\setminus\{x\}$. Therefore
  \begin{equation}
    G^n_{T_\varepsilon^M\setminus\{x\}}(x,x)\le \sum_{k=0}^\infty
    n^{-\gamma k}G^n_{D_x(2^nn^{-\kappa })^c}(x,x)\le\frac {2\log2}\pi n+o(n).
  \end{equation}
  This finishes the proof.
\end{proof}

\textit{Condition~\ref{c:noreturn}.} 
This condition can be easily verified as in the REM case using 
Lemma~\ref{l:poiscloudpropZd} and Proposition~\ref{p:Zdpot}.

\textit{Condition~\ref{c:post}.} The proof of the last condition is 
slightly more complicated, it follows the lines of the proof of 
Lemma~7.4 in \cite{BCM06}. For any $\delta >0$ we have defined event 
$A_n(\delta )$ by (see \eqref{e:defAn})
\begin{equation}
  \exists 1<j_n<\zeta_n: 
  S_n(r_n(j_n))\le t'_n\le S_n(r_n(j_n+1))-\delta t(n),
\end{equation}
where $t(n)/2\le t'_n\le (1+\theta )t(n)$ is a deterministic time sequence.

We want to show that 
$\mathbb P[X_n(t'_n)=U_n(j_n)|A_n(\delta ),\boldsymbol \tau ]>1-\delta$ 
for all $n$ large enough. 
To simplify the notation we define $y=U_n(j_n)$, $u_n=S_n(r_n(j_n))$. 
We further fix $\omega =(\gamma +7)/(1-\alpha )$ (this value has no 
  particular importance, any larger value would work), and 
$\lambda = \gamma + \omega \alpha /2 +1$.

By Condition~\ref{c:shallow} it is possible to choose $\varepsilon $ 
small enough such that the mean time spent in $T^\varepsilon $ before 
$\xi_n$ is smaller than $\delta^3 t(n)$. This implies (using 
  Chebyshev inequality) that 
\begin{multline}
  \mathbb P\big[\exists v_n\in [S_n(r_n(j_n+1))-\delta t(n), S_n(r_n(j_n+1))]:\\
    X_n(s_n)=y|A_n(\delta ),\boldsymbol \tau  \big] \ge 1-c\delta^2.
\end{multline}
It is easy to show that conditionally on the event in the last 
display, the probability that the process $X_n$ leaves disk 
$\mathcal D:=D_{y}(2^{n}n^{-\lambda })$ between times $u_n$ and $v_n$ 
tends to $0$ as $n$ increases. Therefore,
\begin{equation}
  \mathbb P\big[X_n(t')\in \mathcal D\,\forall t'\in[u_n,v_n]
    |A_n(\delta ),\boldsymbol \tau \big]\ge 1-c\delta^2.
\end{equation}
We use $B_n$ to denote the event in the last display. To finish the 
proof of Condition~\ref{c:post} it is sufficient to show that for all 
$n$ large
\begin{equation}
  \mathbb P[X_n(t'_n)=y|
    \boldsymbol \tau , B_n, u_n, y]\ge 1-\delta /2.
\end{equation}

The Markov process 
$(X_n(u_n+s) : s\in[0,v_n-u_n] )$, given 
$\boldsymbol \tau$ , $B_n$, $u_n$, $y$ is equal in law 
to the process $(V(s): s \in [0, v_n-u_n]) $ conditioned on the 
event $\{ T > v_n-u_n \}$ where $V$ and $T$ are constructed as 
follows:

(i) $V$ stays at site $y$ for an exponential, mean 
$\tau_{y}$, amount of time, then 

(ii) with probability $p(n)$, the probability that a random walk 
starting at $y$ escapes $\mathcal D$ 
before returning to site $y$, the process terminates and $T$ is the 
termination time.  With probability $1-p(n)$ the process $V$ performs 
an excursion away from $y$ conditioned not to leave $\mathcal D$.  At 
the end of the excursion it returns to $y$ and step (i) resumes and 
so on. 

The important point is that the number $p(n)$ is of order $1/n$ while 
(recall $y \in T^M_ \varepsilon $) the mean time spent at $y$ per 
visit exceeds $\varepsilon 2^{2n/\alpha}n^{-\gamma/\alpha  }$.  Thus the conditioning 
event has probability bounded below by $C(\varepsilon , \theta)$.  
Hence it will suffice to show that 
$\mathbb P[\bar V(t'_n-u_n  ) \ne y |\boldsymbol\tau]$ tends to zero as 
$n$ tends to infinity $\boldsymbol\tau  $-a.s.~where process 
$(\bar V(u): u \geq 0 )$ is a Markov process that alternates staying 
at site $y$ an exponential amount of time with mean $\tau_{y} $ and 
performing excursions away from $y$ conditioned to stay within 
$\mathcal D$ (again staying at each site a time according to 
  $\boldsymbol\tau$).

We first show that $\boldsymbol \tau$-a.s. for $n$  
sufficiently large, the expected duration of a conditioned excursion 
from $y$ is very small compared to $\tau_y$ uniformly over possible $y \in T^M_\varepsilon$.
It is easy to prove that in $\mathcal D$ 
there are only traps shallower than 
$\varepsilon n^{-\omega }2^{2n/\alpha }n^{-\gamma /\alpha }$. 
Indeed, let $\mathcal B(y)$ be the event
\begin{equation}
  \mathcal B(y)=\Big\{y\in T^M_\varepsilon , \exists x\in \mathcal D, 
    \tau_x\ge\varepsilon  n^{-\omega} {2^{2n/\alpha }}n^{-\gamma 
      /\alpha }\Big\}.
\end{equation}
Then,
\begin{equation}
  \mathbb P[\bigcup_{y\in \mathcal V_n}\mathcal B(y)]\le C 2^{2n} 
  (2^{-2n}n^\gamma )^2 
  n^{\omega \alpha  }2^{2n}n^{-2\lambda } = C n^{-2}
\end{equation}
and the claim follows easily by the Borel-Cantelli lemma.

It was proved, e.g.~in \cite{BCM06}, that 
the expected number of visits to 
$z\in \mathcal D\setminus\{y\}$ during an excursion that does not leave the 
disk is smaller than $2$. The expected duration of the $i$-th excursion, $V_i$, thus satisfies
\begin{equation}
  \mathbb E[V_i|\boldsymbol \tau]\le 2\sum_{z\in \mathcal D\setminus\{y\}}\tau_z
  \le 2\sum_{z\in \mathcal V_n}\tau_z\bbone\{\tau_z\le
    n^{-\omega }\varepsilon 2^{2n/\alpha }n^{-\gamma /\alpha }\}.
\end{equation}
It follows from \eqref{e:tauZ} that 
\begin{equation}
  \mathbb E[\tau_x\bbone\{\tau (x)\le u\}]= O( u^{1-\alpha })\qquad(u\to\infty).
\end{equation}
Therefore, by Chebyshev inequality, 
\begin{equation}
  \mathbb P\Big[\sum_{z\in \mathcal V_n}\tau_z\bbone\{\tau_z\le
      n^{-\omega }\varepsilon 2^{2n/\alpha }n^{-\gamma /\alpha }\}\ge 
    2^{2n/\alpha }n^{-\gamma/\alpha  }n^{-5}\Big]\le C n^{-2}.
\end{equation}
Hence, for a.e.~$\tau $ and $n$ large enough 
\begin{equation}
  \mathbb E[V_i|\boldsymbol \tau ]\le C 2^{2n/\alpha }n^{-\gamma 
    /\alpha }n^{-5}.
\end{equation}

The expected number of excursions of $\bar V$ before time 
$v_n-u_n=O(t(n))$ is bounded by a multiple of $n$, the mean of the total 
time spent by $\bar V$ during the interval 
$[0, v_n-u_n]$ away from $y$ is easily bounded 
by $C2^{2n/\alpha }n^{-\gamma /\alpha }n^{-4}$ for $C$ depending on $\varepsilon $ but 
not on $n$.  

We claim that (for $n$ sufficiently large) for any 
$s \in [0, v_n-u_n]$, $\mathbb P[\bar V(s) \ne y] \leq 2C/n^2$. 
Suppose not.  Then for some $s'$, 
$\mathbb P[Y(s') \ne y] \ge 2C/n^2$. We have that the expected 
total time spent by $\bar V$ away from $y$ in interval 
$[s', s'+2^{2n/\alpha}n^{-\gamma /\alpha }n^{-2}]$ is bounded by 
$C2^{2n/\alpha }n^{-\gamma /\alpha }n^{-4}$, 
so there exists $s'' \in [s', s' +2^{2n/\alpha}n^{-\gamma /\alpha }n^{-2}]$ so that 
$\mathbb P[\bar V(s'') \ne y] \leq C/n^2.$  On the other hand, by 
the Markov property for $\bar V$, if $\nu $ is the time of the 
first jump from $y$,
\begin{equation}
  \mathbb P[\bar V(s'') \ne y] \geq \mathbb P[\bar V(s'') \ne y \cap \{ \nu > s''-s' \}]
  > \frac{1}{2} \mathbb P[\bar V(s') \ne y] \geq C/n^2.
\end{equation}
for $n$ sufficiently large.  This contradiction gives the desired result.

\def\cprime{$'$}
\providecommand{\bysame}{\leavevmode\hbox to3em{\hrulefill}\thinspace}
\providecommand{\MR}{\relax\ifhmode\unskip\space\fi MR }
\providecommand{\MRhref}[2]{%
  \href{http://www.ams.org/mathscinet-getitem?mr=#1}{#2}
}
\providecommand{\href}[2]{#2}

\end{document}